\documentclass[10pt]{article}
\usepackage{amsmath, amssymb}
\usepackage{latexsym}
\usepackage{mathtools}
\usepackage{ulem}
\usepackage{multirow}
\usepackage{graphicx}
\usepackage{color}
\usepackage{url}
\usepackage{framed}
\usepackage{float}
\usepackage{enumitem}
\usepackage{lipsum}
\usepackage{appendix}
\providecommand{\keywords}[1]{\textbf{Keywords:} #1}

\numberwithin{equation}{section}
\DeclareMathAlphabet{\itbf}{OML}{cmm}{b}{it}

\newcommand{\RR}{\mathbb{R}}

\newcommand{\ds}{\displaystyle}

\newcommand{\ri}{\rightarrow}

\newcommand{\bm}{{\itbf m}}

\newcommand{\bx}{{\itbf x}}

\newcommand{\bg}{{\itbf g}}

\newcommand{\bev}{{\itbf e}}
\newcommand{\bw}{{\itbf w}}
\newcommand{\bu}{{\itbf u}}
\newcommand{\bn}{{\itbf n}}
\newcommand{\by}{{\itbf y}}

\newcommand{\bi}{\begin{itemize}}
\newcommand{\ei}{\end{itemize}}

\newcommand{\cu}{{\cal U}}
\newcommand{\cv}{{\cal V}}
\newcommand{\cg}{{\cal G}}

\newcommand{\bH}{{\itbf H}}

\newcommand{\be}{\begin{eqnarray}}
\newcommand{\ee}{\end{eqnarray}}

\newcommand{\ben}{\begin{eqnarray*}}
\newcommand{\een}{\end{eqnarray*}}
\def\ds{\displaystyle}

\newcommand\ov{\overline}

\newtheorem{prop}{Proposition}[section]

\newtheorem{thm}{Theorem}[section]

\newcommand{\bea}{\begin{eqnarray*}}
\newcommand{\eea}{\end{eqnarray*}}
\newcommand{\bean}{\begin{eqnarray}}
\newcommand{\eean}{\end{eqnarray}}
\newcommand{\p}{\partial}
\newcommand{\f}{\frac}

\newcommand{\di}{\mbox{div }}
\newcommand{\aaa}{\mbox{$[$}}
\newcommand{\bbb}{\mbox{$]$}}

\begin{document}

\title{A parallel  sampling algorithm 
for  some nonlinear inverse problems}

\author{  Darko
Volkov \thanks{Department of Mathematical Sciences,
Worcester Polytechnic Institute, Worcester, MA 01609.
}  }

\maketitle

\begin{abstract}
We derive 
a parallel  sampling  algorithm for 
computational inverse problems 
that present an unknown  linear forcing term and a vector of nonlinear parameters 
to be recovered. 
It is assumed that the data is noisy and that the linear part of the problem is ill-posed. 
The vector of nonlinear parameters $\bm $ 
is modeled as a random variable. 
A dilation parameter $\alpha$ is used to scale the regularity of the linear unknown
and is also modeled as a random variable. 
A posterior probability distribution for $(\bm, \alpha)$
is derived following an approach related to 
the maximum likelihood 
regularization parameter selection \cite{galatsanos1992methods}.
A major difference in our approach is that,  
unlike in \cite{galatsanos1992methods}, we do not limit ourselves to  the maximum
likelihood value of $\alpha$.
We then derive a parallel sampling algorithm 
where  we alternate computing proposals in parallel and combining proposals to accept or reject them
as  in  \cite{calderhead2014general}. 
This algorithm is well-suited to problems where proposals are expensive to compute.
We then apply it  to an inverse problem in seismology. 
We show how our results compare favorably to those obtained from 
the Maximum Likelihood (ML),   
the Generalized Cross Validation (GCV), and the Constrained Least Squares (CLS)
algorithms.
\end{abstract}
\keywords{Regularization,  Linear and nonlinear inverse problems, Markov chains, Parallel computing, Elasticity equations in unbounded domains.}

\bigskip

\section{Introduction}
Many physical phenomena are modeled by governing equations 
that depend linearly on some terms and non-linearly on other terms.
For example, the wave equation may depend linearly on a forcing term
and non-linearly on the medium velocity.
This paper is on inverse problems where both a linear part and a nonlinear
part are unknown. Such inverse problems  occur
in passive radar imaging, or in seismology where the source of an earthquake 
has to be determined (the source could be a point, or a fault) and 
a forcing term supported on that source is also unknown.
This inverse problem is then linear in the unknown forcing term and
nonlinear in the location of the source.  
The training phase of neural 
networks is another instance of an inverse problem
 where both
 linear and nonlinear  unknowns occur.
Most neural networks  are based on  linear combinations of basis functions 
depending on a few parameters and 
these parameters have to be determined by training the network on data,
in other words by solving an 
inverse problem that combines linear and nonlinear unknowns
\cite{bishop1995neural}.\\
Let us now formulate the inverse problem studied in this paper. 
 Assume that after discretization a  model 
leads
to the relation
\bean
 \bu = A_{\bm} \bg+ {\cal E}, \label{beq} 
\eean
where $\bg$ in $\RR^p$ is a forcing term, 
$\bm$ in ${\cal B } \subset \RR^q$  is a nonlinear parameter, $A_\bm$ is an $n \times p$ matrix depending continuously on the 
parameter $\bm$,  ${\cal E}$ is  a
Gaussian random variable in $\RR^n$ modeling noise, and $\bu$ in $\RR^n$ is the resulting data for the inverse problem. 
For a fixed $\bm$,
this is the same model that 
Golub et al.  considered in \cite{golub1979generalized}. 
In our case, our goal is to estimate  the nonlinear parameter $\bm$ from $\bu$.
%
We assume that
 the mapping from ${\cal B }$ to $\RR^{n \times p}$,
$\bm \ri A_\bm$, is  known, in other words, a model is known. 
We are interested in the challenging case where  the following difficulties arise simultaneously:
\begin{enumerate}[label=(\roman*)]
\item the size $n \times p$ of the matrix $A_{\bm}$ is such that $n << p$ ({\sl sparse data}),
\item the  singular values of $A_{\bm}$, $s_1 \geq ... \geq s_n>0$ are such that 
   $s_1 >> s_n >0$ ({\sl $A_{\bm} A_{\bm}'$ is invertible but ill-conditioned}),
	\item ${\cal E}$ has zero mean and covariance $\sigma^2 I_n$, but $\sigma$ is unknown.
\end{enumerate}
Due to (i) and (ii),  
$\| A_\bm  \bg - \bu \|$ can be made 
	arbitrarily small for some $\bg$ in $\RR^p$, thus the Variable Projection (VP) functional as defined
	in \cite{golub2003separable} can be minimized
	to numerical zero for all $\bm$ in the search set, and
	the Moore-Penrose 
	inverse can not be used for this problem.
Additionally, a numerical algorithm only based on   
minimizing a  regularized functional may not adequately take into account
the noise, which is amplified by ill-conditioning and nonlinear effects,
and such an algorithm can easily get trapped in local minima.
We thus set up the inverse problem 
consisting of finding $\bm$ from $\bu$ using \eqref{beq}
in a Bayesian framework where we seek to compute
the posterior distribution
of $\bm$. 
%
The inherent advantage of this probabilistic approach is that
related  Markov chains algorithms can avoid being trapped in local minima
by occasionally accepting proposals of lower probability.\\
%
Stuart provided in \cite{stuart2010inverse}
 an extensive survey of probabilistic methods for inverse problems derived from PDE
models and established connections between continuous formulations and their discrete
equivalent. 
Particular examples found in \cite{stuart2010inverse} 
include an inverse problem for a diffusion coefficient, recovering the 
initial condition for the heat equation, and 
determining the permeability of subsurface rock using 
Darcy's law.
In the  application to geophysics that we cover
in section \ref{num sim},
$\bg$ models the slip on a  fault, and $\bm$ is a parameter for 
modeling the piecewise linear geometry of a fault.
The case of interest in this paper is particular 
due to the combination 
of the linear unknown $\bg$ which lies in $\RR^p$, where $p$ is large, 
and the nonlinear unknown $\bm$  which lies in $\RR^q$, with $q<<p$. 

 Tikhonov regularization  may be used instead
of the Moore-Penrose inverse
 to avoid a high norm or a highly oscillatory $\bg$ in
\eqref{beq}.
However, how regular solutions should be is unclear due to (iii).
Accordingly, we  introduce the regularized error functional,
\bean
 \|  A_\bm \bg -  \bu\|^2   + \alpha \| R \bg \|^2,   \label{reg0}
\eean
where 
$R$ is an invertible $p$ by $p$ matrix and $\alpha >0$ is a scaling parameter.
Typical choices for $R$ include the identity matrix and 
 matrices derived from  discretizing   derivative operators.
Without loss of generality, we can 
consider the functional
\bean
 \|  A_\bm \bg -  \bu\|^2   + \alpha \| \bg \|^2  \label{reg},
\eean
in place of\eqref{reg0} by redefining  $A_\bm$ as $A_\bm R^{-1} $.\\
Our solution method will rely on a Bayesian approach.
Assuming that the prior of $\bg$ is also Gaussian, it is well-known that 
the functional
(\ref{reg}) can be related 
to the probability density of $\bu$ knowing $\sigma, \bm$ and $ \alpha$.
However, $\sigma$ is unknown. 
 We will  use the Maximum Likelihood (ML) assumption to eliminate $\sigma$.
As far as we know, this idea was first introduced (for linear problems only) in \cite{galatsanos1992methods},
but unlike  in that reference, we do not eliminate $\alpha$.
We let  $\alpha$  be a random 
variable and thanks to Bayes' theorem we find  
a formula for the  probability density of $(\bm, \alpha)$ knowing $\bu$: this is
 stated in proposition
\ref{mainth}.
In section \ref{par}, 
this formula  is used 
to build a parallel adaptive  sampling algorithm  
to simulate the probability density of $(\bm, \alpha)$ knowing $\bu$.
Finally we show in section  \ref{num sim}  numerical simulations
where this algorithm is applied to a particularly
challenging inverse problem in geophysics.
In this problem a fault geometry described by a nonlinear parameter 
$\bm$ in $\RR^6$  has to be reconstructed from surface displacement data modeled by 
the vector $\bu$ in $\RR^n$. The data is produced by a large slip field ${\cal G}$
modeled by a vector $\bg$  in $\RR^p$.
$\bu$ depends linearly on $\bg$ and this dependance can be expressed by a matrix
$A_{\bm}$. In this simulation, the matrix $A_{\bm}$ is full since it is derived
from convolution by a Green function. In addition, the entries of $A_{\bm}$
are particularly expensive to compute, which is a hallmark of problems involving
 half space elasticity and
this application features all the difficulties (i), (ii), and (iii) listed above.
This makes the VP functional method unsuitable,
and it also renders
classical minimization methods such as   GCV, CLS, and ML
much less accurate than our proposed method, 
as shown in section
\ref{comparisons}.



\section{Solution method}

\subsection{The linear part of the inverse problem and selection methods for $\alpha$}
There is a vast amount of literature on methods for selecting an adequate value  for the regularization parameter $\alpha$, assuming that the nonlinear parameter $\bm$ is fixed.
An account of most commonly used methods, together with error analysis,
 can be found in \cite{vogel2002computational}.
In this paper we review  three such methods, that we later compare 
to our own algorithm.
Throughout the rest of this paper,
the Euclidean norm will be denoted by $\| . \|$ and the transpose of a matrix $M$
will be denoted by $M'$.

\subsubsection{Generalized cross validation (GCV)}
The GCV method was first introduced and analyzed in 
\cite{golub1979generalized}. 
The parameter $\alpha$ is selected by minimizing
\bean
\f{\|  (I_n - A_\bm(A_\bm'A_\bm + \alpha I_p)^{-1} A_\bm') \bu \|^2}{ \mbox{tr }(I_n - 
A_\bm(A_\bm'A_\bm + \alpha I_p)^{-1} A_\bm')^2}, \label{GCV}
\eean
where tr is the trace.
Let $\alpha_{GCV}$ be the value of $\alpha$ which minimizes (\ref{GCV}).
Golub et al. proved in \cite{golub1979generalized} that
as a function of $\alpha$,
the expected value of $ \|  A_\bm \bg -  A_\bm(A_\bm'A_\bm + \alpha I_p)^{-1} A_\bm'\bu \|^2$
which can be thought of as an indicator of fidelity of the pseudo-solution 
$(A_\bm'A_\bm + \alpha R'R)^{-1} A_\bm'\bu $ 
is approximately minimized at $\alpha=\alpha_{GCV}$ as $n \ri \infty$.
Although the GCV method enjoys this remarkable asymptotic property
and does not require knowing  $\sigma^2$, many authors have noted that 
determining the minimum of (\ref{GCV}) in practice can be costly and inaccurate as  in
 practical situations
the quantity in (\ref{GCV})  is  flat near its minimum for a wide range of values of 
$\alpha$ \cite{thompson1989cautionary, varah1983pitfalls}.

\subsubsection{Constrained least square (CLS)}
This method, also called  
the discrepancy principle \cite{morozov1966solution, vogel2002computational},
advocates choosing a value for $\alpha$ such that 
\bean 
\|  \bu -  A_\bm(A_\bm'A_\bm + \alpha I_p)^{-1} A_\bm'\bu \|^2 = n \sigma^2.\label{CLS}
\eean
Clearly, applying this method requires knowing 
$\sigma^2$ or at least some reasonable approximation of its value.
Even if  $\sigma^2$ is known, this method leads to solutions that are in general overly
smooth, \cite{galatsanos1992methods, vogel2002computational}.

\subsubsection{Maximum likelihood (ML)}\label{ML method}
To the best of our knowledge the ML  method was first proposed in \cite{galatsanos1992methods}.
It relies on 
the fundamental assumption that the prior of $\alpha^{\f12} \bg$ is also  normal
with zero mean and covariance $\sigma^2 I_p$.
The likelihood of the minimizer of (\ref{reg}) knowing $\sigma$ and 
$\alpha$ is then maximized for all $\alpha >0$, $\sigma >0$.
Galatsanos and Katsaggelos showed
in \cite{galatsanos1992methods} 
that equivalently the expression
\bean
\f{\bu'(I_n - A_\bm(A_\bm'A_\bm + \alpha I_p)^{-1} A_\bm') \bu }{(\det (I_n - A_\bm(A_\bm'A_\bm + \alpha I_p)^{-1} A_\bm' ))^{1/n}}    \, , \label{MLratio}
\eean
has to be minimized for all $\alpha >0$, which 
does not require knowing $\sigma^2$.

\subsection{A  sampling algorithm for computing the posterior of  $(\bm, \alpha)$}
\label{new Bayesian approach}
In this section 
we  generalize the ML method  recalled in   section \ref{ML method} to problems
depending nonlinearly on the random variable $\bm$, 
while refraining from only
 retaining the maximum likelihood 
value of $\alpha$.

	\begin{prop} \label{mainth}
Assume that
\begin{enumerate}[label=H\arabic*. , wide=0.5em,  leftmargin=*]
  \item $\bu$, $\bg$, $\bm$ and $\alpha$ are random variables in 
	   $\RR^n, \RR^p, {\cal B} \subset \RR^q, (0, \infty)$, respectively,
		\label{first assump}
  \item $(\bm, \alpha)$ has a known prior distribution denoted by
	   $\rho_{pr} (\bm, \alpha)$, \label{bm, alph assump}
  \item $A_{\bm}$ is an $n$ by $p$ matrix which depends continuously on $\bm$,
	\item ${\cal E}$ is an $n$ dimensional normal random variable
	with
	 zero mean  and covariance $\sigma^2 I_n$, \label{cal E assump}
\item relation (\ref{beq}) holds, \label{1-1 assump}
	\item the prior of $\alpha^{\f12} \bg$ is  a normal random variable 
	with zero mean  and covariance $\sigma^2 I_p$. \label{ML assump}
\end{enumerate}
Let $\rho(\bu | \sigma, \bm, \alpha)$ be the conditional probability density of 
$\bu$ knowing $\sigma, \bm, \alpha$. 
As a function of $\sigma>0$, $\rho(\bu | \sigma, \bm, \alpha)$ achieves a unique maximum at
\bean
\sigma_{max}^2 = \f{1}{n} ( \alpha \| \bg_{min} \|^2  +\| \bu - A_\bm \bg_{min} \|^2),
\label{sigma}
\eean
where $\bg_{min} = (A_\bm' A_\bm + \alpha I_p)^{-1} A_\bm' \bu$
 is the minimizer of (\ref{reg}).
Fixing $\sigma=\sigma_{max}$, the probability density of $(\bm, \alpha)$ knowing $\bu \neq 0$ is then
given, up to a multiplicative constant, by the formula
\bean
\rho ( \bm ,\alpha| \bu) \propto  \det (\alpha^{-1} A_\bm'A_\bm + I_p)^{-\f12} ( \alpha \| \bg_{min} \|^2  +\| \bu - A_\bm \bg_{min} \|^2          )^{-\f{n}{2}} \rho_{pr} (\bm, \alpha ) .  \label{final}
\eean
\end{prop}
Note that  Galatsanos and Katsaggelos give a proof in \cite{galatsanos1992methods} 
  of a closely related result, but  our context is somehow different. 
	Our matrix  
	 $A_\bm$ is rectangular and depends continuously and non-linearly on the random variable $\bm$.
	In the present case, a similar argument can be carried out
	noting that thanks to assumption \ref{ML assump},
	the marginal probability density of $\bu$ 
	knowing $\sigma, \bm, \alpha$ can be computed and equals,
	using the minimizer $\bg_{min} = (A_\bm' A_\bm + \alpha I_p)^{-1} A_\bm' \bu$
	of   (\ref{reg}),
	$$
	(\f{1}{2 \pi \sigma^2})^{\f{n}{2}}  \alpha^{\f{p}{2}}
\exp (- \f{\alpha}{2 \sigma^2} \| \bg_{min} \|^2  - \f{1}{2 \sigma^2} \| \bu - A_\bm \bg_{min} \|^2          ) 
(\det ( A_\bm' A_\bm + \alpha I_p))^{-\f12},
	$$
	which is in turn maximized for $\sigma$ in $(0, \infty)$.
	This yields formula \eqref{sigma}. Formula \eqref{final} follows from there.
	
	Note that the fundamental assumption \ref{ML assump} was  introduced in \cite{galatsanos1992methods} and 
	can be thought of as  a way of  restoring a balance between 
	reconstruction fidelity (first term in (\ref{reg}))
	and regularity requirements (second term in (\ref{reg})).
We note that the right hand side of (\ref{final})
can be related to the ML ratio \eqref{MLratio}, first because
as $\bg_{min}= (A_\bm' A_\bm + \alpha I_p)^{-1} A_\bm' \bu$,
a simple calculation will show that
\bea
\| \bu - A_\bm \bg_{min} \|^2  + \alpha \| \bg_{min} \|^2 
= \bu'  (I_n - A_\bm (A_\bm'A_\bm + \alpha I_p)^{-1} A_\bm ') \bu,
\eea
and second because
of the identity
\bean
(\det(\alpha^{-1} A_\bm' A_\bm + I_p))^{-1} 
= \det (I_n  - A_\bm (A_\bm'A_\bm + \alpha I_p )^{-1} A_\bm'). \label{dets}
\eean
which is shown in Appendix \ref{det proof}. 
Next it is important
to note that, as explained in  Appendix \ref{det proof},
$$  I_n - A_\bm( A_\bm' A_\bm + \alpha I_p )^{-1}A_\bm' 
=(I_n + \alpha^{-1} A_\bm A_\bm' )^{-1}. $$
Since $n << p$, this is particularly helpful in numerical calculations,
and it allows to use the formula 
$\bu'(I_n + \alpha^{-1} A_\bm A_\bm')^{-1}\bu$ to compute 
$\| \bu - A_\bm \bg_{min} \|^2  + \alpha \| \bg_{min} \|^2 $.
Introducing  $\bg_{min}$ is  helpful, however, for comparing the solution method 
developed in this 
paper to classical variational regularization methods see, section 1.3 in \cite{vogel2002computational}.



\subsection{Single processor algorithm} \label{single}
Based on \eqref{final}, we define the non-normalized distribution
\bean
{\cal R} ( \bm, \alpha ) =\det (\alpha^{-1} A_\bm'A_\bm + I_p)^{-\f12} ( \alpha \| \bg_{min} \|^2  +\| \bu - A_\bm \bg_{min} \|^2          )^{-\f{n}{2}} \rho_{pr} (\bm, \alpha ). \label{non}
\eean
We use the standard notations
$ {\cal N} (\mu ,   \textsf{$\Sigma$} )$ for a normal distribution with mean $\mu$ and covariance
 \textsf{$\Sigma$}, $ {\cal }U (0 ,1)$ for a uniform distribution in the interval $(0,1)$.
The 
algorithm starts from a point 
$(\bm_1, \alpha_1)$ in $\RR^{q+1}$ such that 
$\rho_{pr} (\bm_1, \alpha_1) >0$ and an initial covariance matrix $\Sigma_0$ obtained 
from the prior distribution. A good choice of the initial point 
$(\bm_1, \alpha_1)$  may have a strong impact on how many sampling 
steps are necessary. 
A poor choice may result in  a very long "burn in" phase where the random walk
is lost in a low probability region.
How to find a good starting point $(\bm_1, \alpha_1)$ depends greatly on the application,
so we will discuss that issue in a later section where we cover a specific example.
Our basic single processor algorithm follows the well established 
\textsl{adaptive} MCMC propose/accept/reject algorithm
\cite{roberts2009examples}.
Let $\beta_j$ be a decreasing sequence in $(0,1)$ which converges to 0.
This sequence is used to weigh a convex combination
between the initial covariance  $\Sigma_0$ and 
the  covariance $\Sigma$ learned from sampling.
The updating of $\Sigma$ need not occur at every step.
Let $N$ be the total number of steps and $N'$
the number of steps between updates of $\Sigma$.
We require that $1<N' <N$.
Finally, the covariance for the proposals is adjusted by a
factor of $(2.38)^2 (q+1)^{-1}$
as recommended in \cite{roberts2009examples, roberts2001optimal}.
It was shown in   \cite{roberts2001optimal} that this scaling leads 
to an optimal acceptance rate.


\begin{framed}
\textbf{Single processor sampling algorithm}
\begin{enumerate}
\item Start from a point $(\bm_1, \alpha_1)$ in $\RR^{q+1}$ and set
  $\Sigma=\Sigma_0$.
\item for $j=2$ to  $N$ do:

\begin{enumerate}
\item if $ j $ is a multiple of $N'$ update 
   the covariance $\Sigma$  by using the points $(\bm_k, \alpha_k), 1 \leq k \leq j-1$,
\item draw  $(\bm^*, \alpha^*)$ 
from $  (\bm_{j-1},  \alpha_{j-1}) +  (1- \beta_j){\cal N} ( 0 , (2.38)^2 (q+1)^{-1} \textsf{$\Sigma $}) 
+ \beta_j {\cal N} ( 0 ,   (2.38)^2 (q+1)^{-1} \textsf{ $\Sigma_0$})$,
\item compute
   $ { \cal R}(\bm^*, \alpha^*) $,
	\item draw $u$ from ${\cal U} (0,1)$,
	\item if $u < \f{{ \cal R}(\bm^*, \alpha^*) }{{ \cal R}(\bm_{j-1},  \alpha_{j-1}) }$ set
	  $(\bm_{j}, \alpha_{j}) =  (\bm^*, \alpha^*)$, else set 
		$( \bm_{j}, \alpha_{j}) =  (\bm_{j-1},  \alpha_{j-1}) $.
\end{enumerate}


\end {enumerate}
\end{framed}

\subsection{Parallel algorithm}  \label{par}
Let  $N_{par}$ be the number of  processing units. 
A straightforward way of taking advantage of multiple processors 
is to generate  $N_{par}$ separate chains of samples using 
the single processor algorithm 
described in  section \ref{single} and then concatenate them.
However, computations can be greatly accelerated by analyzing  the proposals produced by the chains
in aggregate  \cite{calderhead2014general, jacob2011using}.
While in section \ref{single} $(  \bm_j, \alpha_j, ) $ was a $q+1$ dimensional vector,
here we set $\textsf{M}_{j}$ to be  a $q+1$ by $N_{par}$ matrix where
the $k$-th column  will be denoted by
$\textsf{M}_{j}(k)$ and is a sample of the  random variable $(\bm, \alpha)$, $k=1, ..., N_{par}$.
Next, if $j \geq 2$, we assemble an 
$N_{par} +1$ by $N_{par} +1$  transition matrix $T$
from the computed non-normalized densities
$ {\cal R} (\textsf{M}_{j-1}(N_{par}))$ and
$ { \cal R}(\textsf{M}^*(k)) $, $k=1, ..., N_{par}$, 
where $\textsf{M}^*$ is the proposal.
Let $\bw$ be the vector in $\RR^{N_{par} + 1}$
with coordinates
\bea
\bw= \big( {\cal R} (\textsf{M}_{j-1}(N_{par})), {\cal R}(\textsf{M}^* (1)), ... ,
 {\cal R}(\textsf{M}^* (N_{par})) \big).
\eea 
The entries of the transition matrix $T$ are given by the following formula  \cite{calderhead2014general},
\bean
T_{k, l} = \ds \left\{ 
\begin{array}{l} \ds
\f{1}{N_{par}} \min \{  1, \f{w_l}{w_k}  \},  \mbox{ if } k \neq l , \\
\ds 1 - \sum_{ 1 \leq  l \leq N_{par}+1, l \neq k}  T_{k,l},  \mbox{ if } k = l .
\end{array}
\right.
\label{trans}
\eean
Note that for $k=1,...,N_{par}+1 $
the row $T_{k,1}, ..., T_{k,N_{par}+1}$
defines a discrete probability distribution on 
$\{ 1, ..., N_{par}+1\}$.
\begin{framed}
\textbf{Parallel sampling algorithm}
\begin{enumerate}

\item Start from a point $(\bm_1, \alpha_1)$ in $\RR^{q+1}$ and set
  $\Sigma=\Sigma_0$. Set the $N_{par}$ columns of $\textsf{M}_1$
	to be equal to $(\bm_1, \alpha_1)$.
\item for $j=2$ to  $N$ do:

\begin{enumerate}
\item if $ j $ is a multiple of $N'$ update 
   the covariance $\Sigma$  by using the points $\textsf{M}_{k}(l), 1 \leq k \leq j-1,
	1 \leq l \leq N_{par}$,

\item  for $k=1$ to $N_{par}$, draw  the proposals $\textsf{M}^*(k)$ 
from $ \textsf{M}_{j-1} (k) +  (1- \beta_j){\cal N} ( 0 ,  (2.38)^2 (q+1)^{-1}\textsf{$\Sigma $}) 
+ \beta_j {\cal N} ( 0 ,  (2.38)^2 (q+1)^{-1}\textsf{$\Sigma_0 $})$,
\item compute  \textbf{in parallel}
   $ { \cal R}(\textsf{M}^*(k)) $, $k=1, ..., N_{par}$,
		\item assemble the $N_{par} +1$ by $N_{par} +1$  transition matrix $T$ as indicated above,
	\item for $k=2,...,N_{par}+1 $ draw an integer $p$
	  in $\{ 1, ..., N_{par}+1\}$ using the  probability distribution  $T_{k,1}, ..., T_{k,N_{par}+1}$;
		 if $p=1$ set $\textsf{M}_{j}(k-1) = \textsf{M}_{j-1}(N_{par})$ (reject), otherwise
		 set $\textsf{M}_{j} (k-1) = \textsf{M}^* (p-1)$ (accept).
\end{enumerate}

\end {enumerate}
\end{framed}
Note that this parallel algorithm is especially well suited to applications 
where computing the non-normalized density $ { \cal R}(\bm, \alpha) $ is  expensive.
In that case, even the naive parallel algorithm
where $N_{par} $ separate chains are computed in parallel will be about
$N_{par} $ times more efficient than the  single processor algorithm.
The parallel algorithm presented in this section is in fact even more efficient
due to superior mixing properties and sampling performance:
the performance is not overly sensitive to tuning of proposal parameters
\cite{calderhead2014general}.

\section{Application to the fault inverse problem in seismology and 
numerical simulations} \label{num sim}
Using  standard rectangular coordinates, let $\bx= (x_1, x_2, x_3)$ denote elements of 
 $\RR^{3}$.
We define $\RR^{3-}$ to be the open half space $x_3<0$.
We use the equations of linear elasticity with
Lam\'e constants $\lambda$ and $\mu$ 
such that $\lambda > 0$ and $\lambda + \mu> 0$.
For a vector field $\cv = (\cv_1, \cv_2, \cv_3)$,
the  stress vector in the 
direction $\bev \in \mathbb R^3 $ will be denoted by
\bea
T_\bev \cv = \sum_{j=1}^3
\left(  \lambda \, \di \cv \, \delta_{ij} + \mu \, (\p_i \cv_j + \p_j \cv_i ) \right)
 e_j.
\eea
Let $\Gamma$ be a Lipschitz open surface which is strictly included in $\RR^{3-}$, 
with  normal
vector $\bn$.
We define the jump  $\aaa \cv \bbb$ of the  vector field 
$\cv$ across $\Gamma$ 
to be
$$
\aaa \cv \bbb (\bx) = \lim_{h \ri 0^+} \cv  (\bx + h \bn) - \cv  (\bx - h \bn),
$$
for $\bx$ in $\Gamma$, if this limit exists.
Let $\cu$ be the displacement field solving
\bean
\mu \Delta \cu+ (\lambda+\mu) \nabla \di \cu= 0  \mbox{ in } \RR^{3-} 
\setminus \Gamma \label{uj1}, \\
 T_{\bev_3} \cu =0 \mbox{ on the surface } x_3=0 \label{uj2}, \\
 T_{\bn} \cu  \mbox{ is continuous across } \Gamma \label{uj3}, \\
 \aaa  \cu\bbb =\cg \mbox{ is a given jump across } \Gamma , \label{uj3andhalf}\\
\cu (\bx) = O(\f{1}{|\bx|^2}),  \nabla \cu (\bx) = O(\f{1}{|\bx|^3}), \mbox{ uniformly as } 
|\bx| \rightarrow \infty,
\label{uj4}
\eean
where 
$\bev_3$ is the vector $(0,0,1)$.

Let  $D$ be a bounded domain in $\RR^{3-}$ 
with   Lipschitz boundary $\p D$ containing $\Gamma$. 
Let 
  $\widetilde{H}^{\f12}(\Gamma)^2$ 
  be the space  of restrictions to $\Gamma$ of tangential fields  in
   $H^{\f12}(\p D)^2$ supported in $\ov{\Gamma}$.
In \cite{volkov2017reconstruction}, we defined the functional space 
${\bf S}$ of vector fields $\cv $ defined in $\RR^{3-}\setminus{\ov{\Gamma}}$
such that $\nabla \cv $ and $\ds \f{\cv}{(1+ r^2)^{\f12}}$ 
are in $L^2(\RR^{3-}\setminus{\ov{\Gamma}})$ and we proved the following
	existence and uniqueness result.
\begin{thm}\label{direct exist unique}
 Let $\cg$ be in $\widetilde{H}^{\f12}(\Gamma)^2$.
The problem
(\ref{uj1}-\ref{uj3andhalf})
 has a unique solution in ${\bf S}$.
In addition, the solution $\cu$ satisfies the decay conditions
(\ref{uj4}).
\end{thm}
Can  both $\cg$ and $\Gamma$  be determined from the data
 $\cu$ given only on the plane $x_3=0$?
The following Theorem shown in \cite{volkov2017reconstruction}
asserts that this is possible if the data is known on a relatively open
set  of the plane $x_3=0$.
\begin{thm}{\label{uniq1}}
Let $\Gamma_1$ and $\Gamma_2$ be two 
connected open 
surfaces that are  unions of two polygons.
For $i$ in $\{ 1,  2\}$, assume that $\cu^i$ solves  (\ref{uj1}-\ref{uj4}) for 
$\Gamma_i$ in place of $\Gamma$ and $\cg^i$, a tangential field in 
$\tilde{H}^{\f12}(\Gamma_i)^2$, in place of $\cg$.
Assume that $\cg^i$ has full support in $\Gamma_i$, that is, 
$\mbox{supp } \cg_i = \ov{\Gamma_i}$.
Let $V$ be a non empty open subset in $\{x_3 =0\}$.
If $\cu^1$ and $\cu^2$ are equal in $V$, then
$\Gamma_1 =  \Gamma_2$
and $\cg^1=\cg^2$.
\end{thm}
Theorems  \ref{direct exist unique} and \ref{uniq1} 
were  proved in \cite{volkov2017reconstruction} 
for media with constant Lam\'e coefficients. Later,
in \cite{aspri2020analysis}, 
the direct problem (\ref{uj1}-\ref{uj3andhalf})
 was analyzed under weaker regularity conditions 
for $\cu$ and $\cg$. In \cite{aspri2020arxiv},
the direct problem (\ref{uj1}-\ref{uj3andhalf})
was proved to be uniquely solvable in case of
piecewise Lipschitz coefficients and general elasticity tensors.
Both \cite{aspri2020arxiv}  and \cite{aspri2020analysis}
include a proof of  uniqueness  for the  fault inverse problem under appropriate assumptions.
In our case, the  solution $\cu$ to problem (\ref{uj1}-\ref{uj3andhalf})  can also be written out 
 as the convolution on $\Gamma$
\bean
  \cu(\bx) = \int_\Gamma \bH(\bx, \by, \bn) \cg(\by) \, d \sigma (\by) \label{int formula},
\eean
where $\bH$ is  the  Green's tensor  associated to the system (\ref{uj1}-\ref{uj4}), and
  $\bn$ is the normal to $\Gamma$. 
The practical determination of this adequate half space Green's tensor $\bH$ was first studied
in  \cite{Okada} and later, more rigorously, in \cite{volkov2009double}.
Due to formula (\ref{int formula}) 
we can define a continuous mapping ${\cal M}$ from tangential fields
$\cg $ in $H^{1}_0(\Gamma)^2$ to surface displacement fields $\cu(x_1, x_2, 0)$
in $L^2(V) $.
Theorem \ref{uniq1} asserts that 
this mapping is injective, so an inverse operator can be defined. 
It is well known, however, that such an operator ${\cal M}$ is compact, therefore
its inverse is unbounded. 

\begin{figure}[htbp]
   \centering
        \includegraphics[scale=.4]{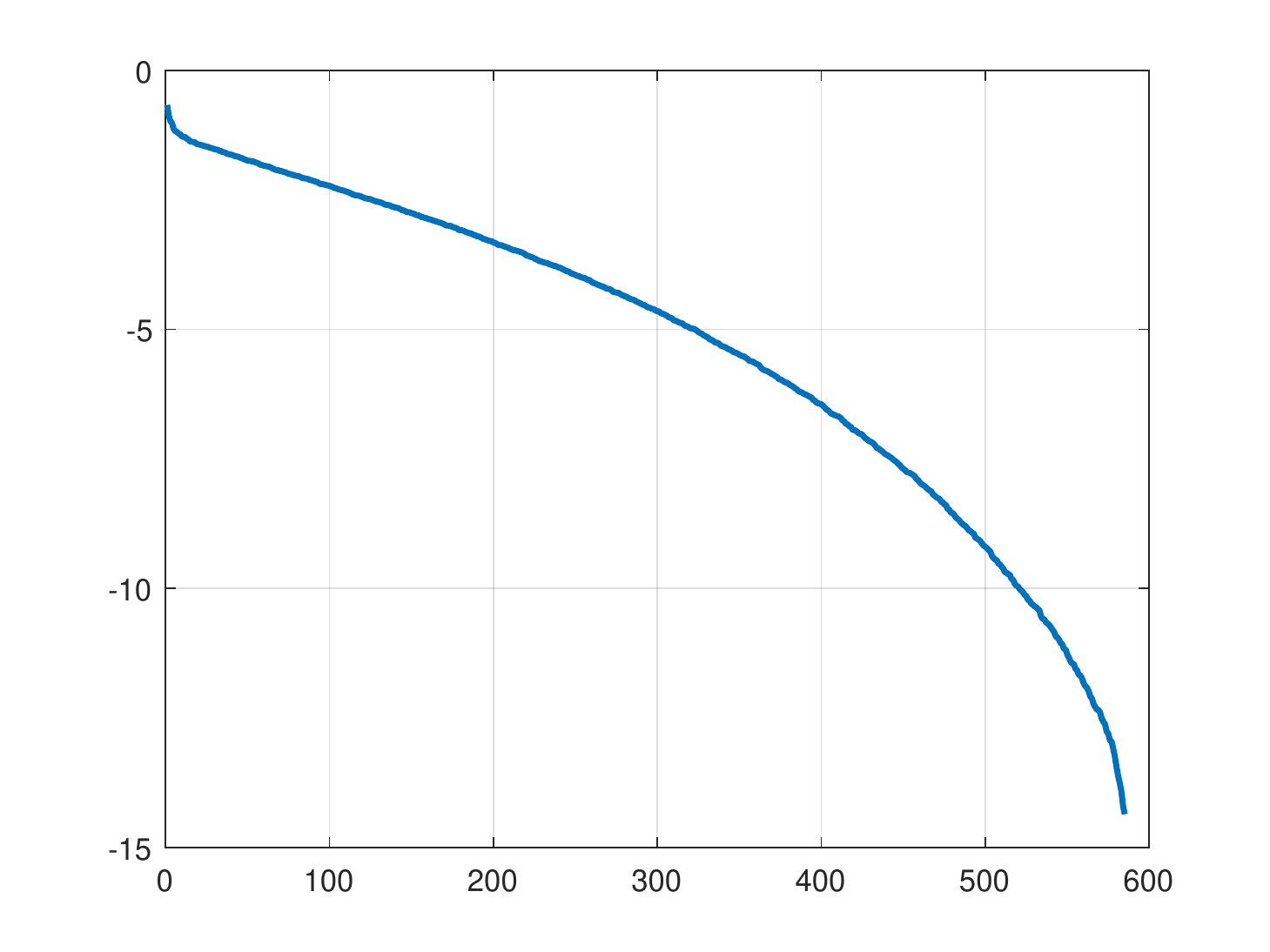} 
    \caption{The decimal logarithm of the singular values of $A_\bm$ with  
		$ n \times p = 585\times 10201$, 
		and   $\bm$ given by  \eqref{values}.  }
    \label{svd}
\end{figure}

We assume here that  $n$ is a multiple of three 
since three-dimensional 
 displacements are measured.  Let 
 $O_j, j=1,..,\f{n}{3}$ be the points on the plane $x_3=0$ where
measurements are collected, giving rise 
 to a data vector
$\bu$ in $\RR^n$.
In our numerical simulations, we assume that 
$\Gamma$ is made up of two contiguous quadrilaterals,
and that  $\Gamma$ is the image by a  piecewise affine function 
of the square $S= [-100, 200] \times [-100, 200]$ in the $x_1  x_2$ plane. 
Applying a change of variables in the integral given by \eqref{int formula}, the field $\cg$ 
may be assumed to be defined on $S$ and the integral itself becomes an integral on $S$. 
We use a regular $m \times m$ grid of points on $S$,
thus we
set $p=m^2$.
The points on this grid are then labeled $Q_k, k=1,..,p$.
Let $\tilde{Q}_k$  be the points on $\Gamma$ with same 
$x_1$ and $x_2$ coordinates as $Q_k$, $ k=1,..,p$:
 $\bg$ in $\RR^p$ is used to approximate
$\cg$ at these points. The Green tensor $\bH$ is evaluated at
 $(\bx,\by)=(O_j,\tilde{Q}_k),  j=1,..,\f{n}{3}, k=1,..,p$, and the integral in 
\eqref{int formula} is approximated by quadrature.
Since $\Gamma$ is the image of $S$ by a piecewise affine function, assuming that 
$\Gamma$  is made up of two contiguous quadrilaterals, it
can be defined by a parameter $\bm$ in $\RR^6$.
We then write
the discrete equivalent of the right hand side of  formula \eqref{int formula} 
as the matrix-vector product $A_\bm \bg$, where
$\bg$ in $\RR^p$ is the discrete analog of $\cg$ and multiplying by the matrix 
$A_\bm $ is the discrete
analog of applying the convolution product of $\cg$ against $\bH$  over $\Gamma$.
Taking into account measurement errors, we arrive 
at the formulation \eqref{beq}. 
In the simulations shown in this paper
the size of the matrix $A_\bm$
is   $ n \times p \sim 500 \times 2500$. 
The singular values of $A_\bm$ decay fast
(this is due to the fast decay of  the singular values of the compact operator
${\cal M}$, see \cite{little1984eigenvalues}), so even 
choosing a coarser grid on $\Gamma$ which would make $p \leq n$ 
would still result in an ill-conditioned matrix $A_\bm' A_\bm$.
In Figure \ref{svd} we plot the singular values of $A_\bm$ for the particular
value of $\bm$ used to generate the direct data for the inverse problem
used for illustration in the next section. 
A similar fast decay of these singular values is observed for all $\bm$
in the support of its prior.

Another practical aspect of the matrix $A_m$ is that 
it  is full (as  is usually the case in problems derived from integral operators) and its entries are  expensive to compute (this is due to the nature of the half space elastic Green tensor)
\cite{volkov2009double}.  However, great gains can be achieved by applying array 
operations thus
taking advantage of multithreading.
The matrix $R$ used to regularize $\bg$ as in \eqref{reg0} 
is based on derivatives in the $x_1$ and $x_2$ directions.
We can find  two $p$ by $p$ permutation matrices $P_D$ and $P_E$ 
satisfying $D = P_D^{-1} M P_D  $, $E = P_E^{-1} M P_E  $, 
  where $M$ is block-diagonal with each
block of size $m \times m$ given by 
$$\begin{bmatrix}
1      & -1     & 0     &      0  &\cdots         & 0  \\
0      & 1      & -1    & 0       &    \cdots     & 0\\
\vdots & \ddots & \ddots & \ddots &   \ddots      & \vdots    \\
\vdots & \ddots & \ddots & \ddots &   \ddots      & 0    \\
 0     & \cdots & \cdots & 0 &   1           & -1      \\
 0     & \cdots & \cdots & \cdots &   0           & 1
\end{bmatrix},$$ 
 such  that $D$ is the discrete analog of a derivative operator in the $x_1$ direction
and $E$ is the discrete analog of the  derivative operator in the $x_2$ direction.
Note that $\| D \| \leq 2 $, $\| E\| \leq 2 $, 
and $\| D^{-1} \| \leq p $, $\| E^{-1}\| \leq p $, in matrix 2 norm.  
We then define $R$ to be a  matrix such that  $R'R=D' D  + E' E  $. 
Evaluating $R$ is unnecessary since 
 only $R'R$ is used in  computations.
For efficiency, it is advantageous not to evaluate the matrix product $A_\bm R^{-1}$ 
to reduce \eqref{reg0} to \eqref{reg}.
Instead, we evaluate and store $R'R$ 
and we use an iterative solver 
to evaluate $\bg_{min} =  (A_\bm' A_\bm + \alpha R'R)^{-1} A_\bm' \bu$.
We coded the
 function $ \bg \ri (A_\bm' A_\bm + \alpha R'R) \bg$ without evaluating
the matrix product $A_\bm' A_\bm$.

\subsection{Construction of the data}
\begin{figure}[htbp]
    \centering
      \includegraphics[scale=.5]{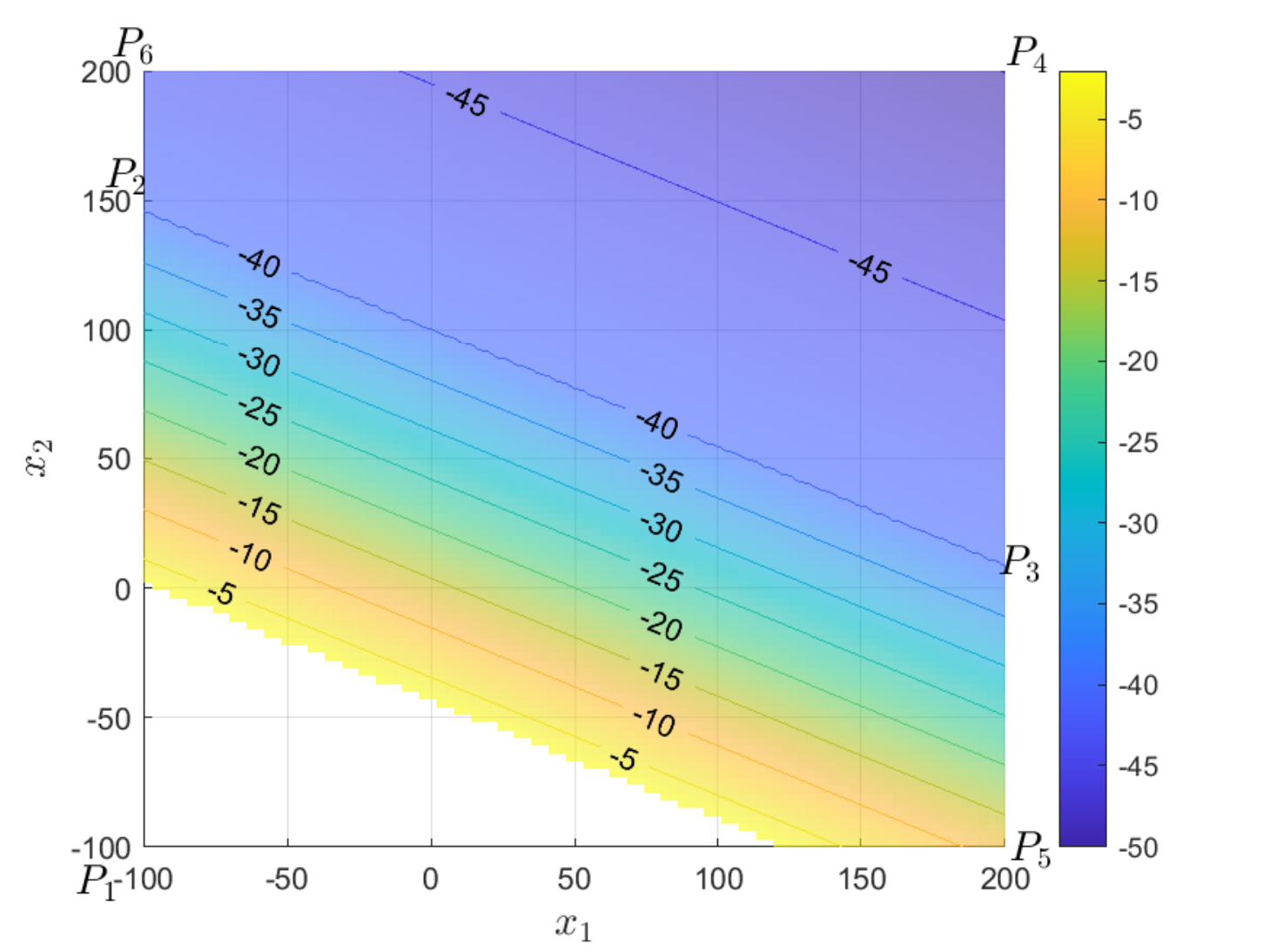}
      \includegraphics[scale=.5]{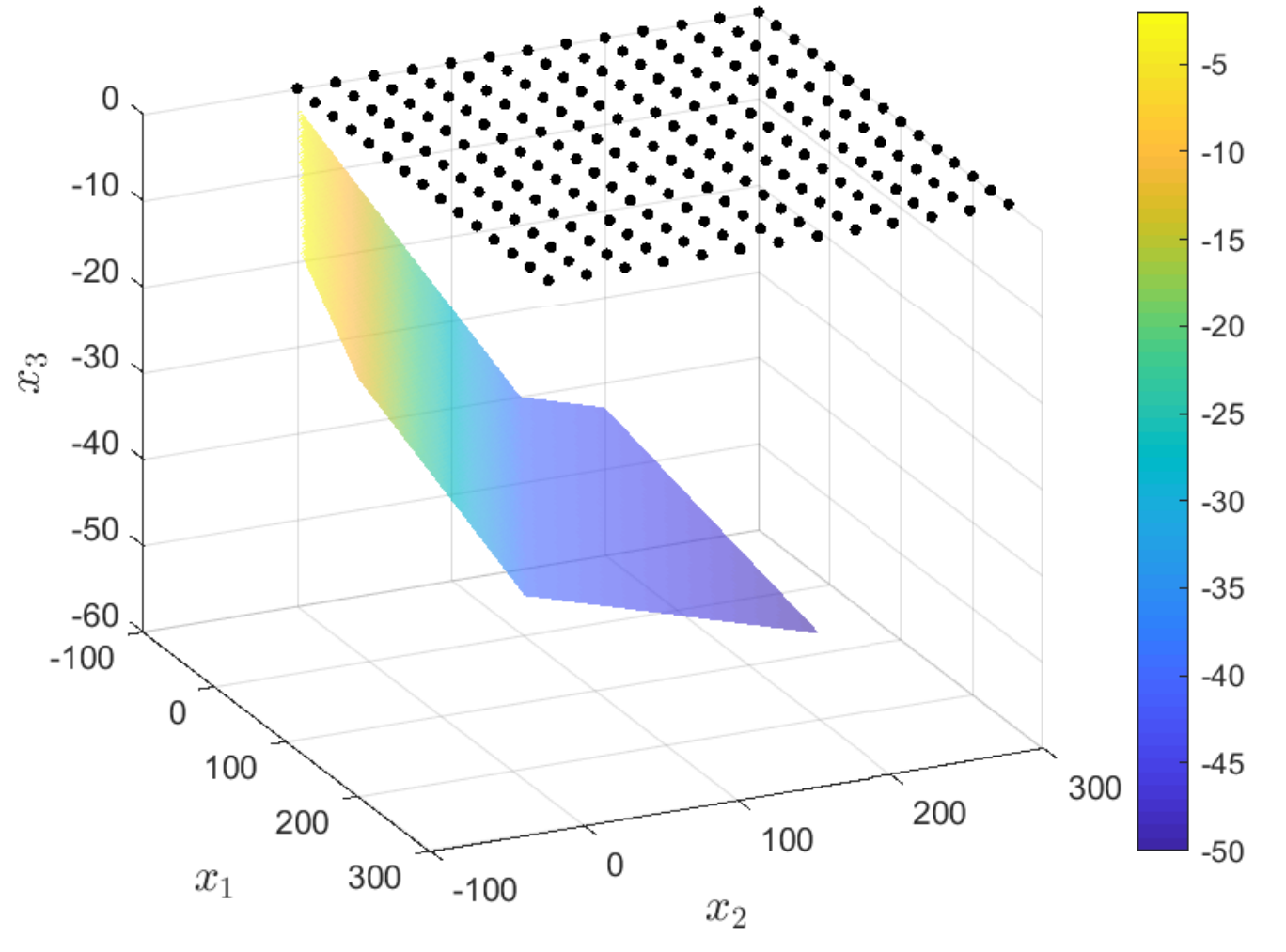}
    \caption{The piecewise planar  surface  $\Gamma$.
		Depths are indicated by  color bars.  
		Top graph: view of $\Gamma$ from above with the six points
		$P_1, P_2, P_3, P_4,  P_5, P_6$ and contour lines of same depth.
		Bottom graph: a three dimensional rendition of $\Gamma$.
	The measurement points are on the surface $x_3=0$ and
		are indicated by black dots.}
    \label{fault}
\end{figure}

\begin{figure}
   \centering
        \includegraphics[scale=.5]{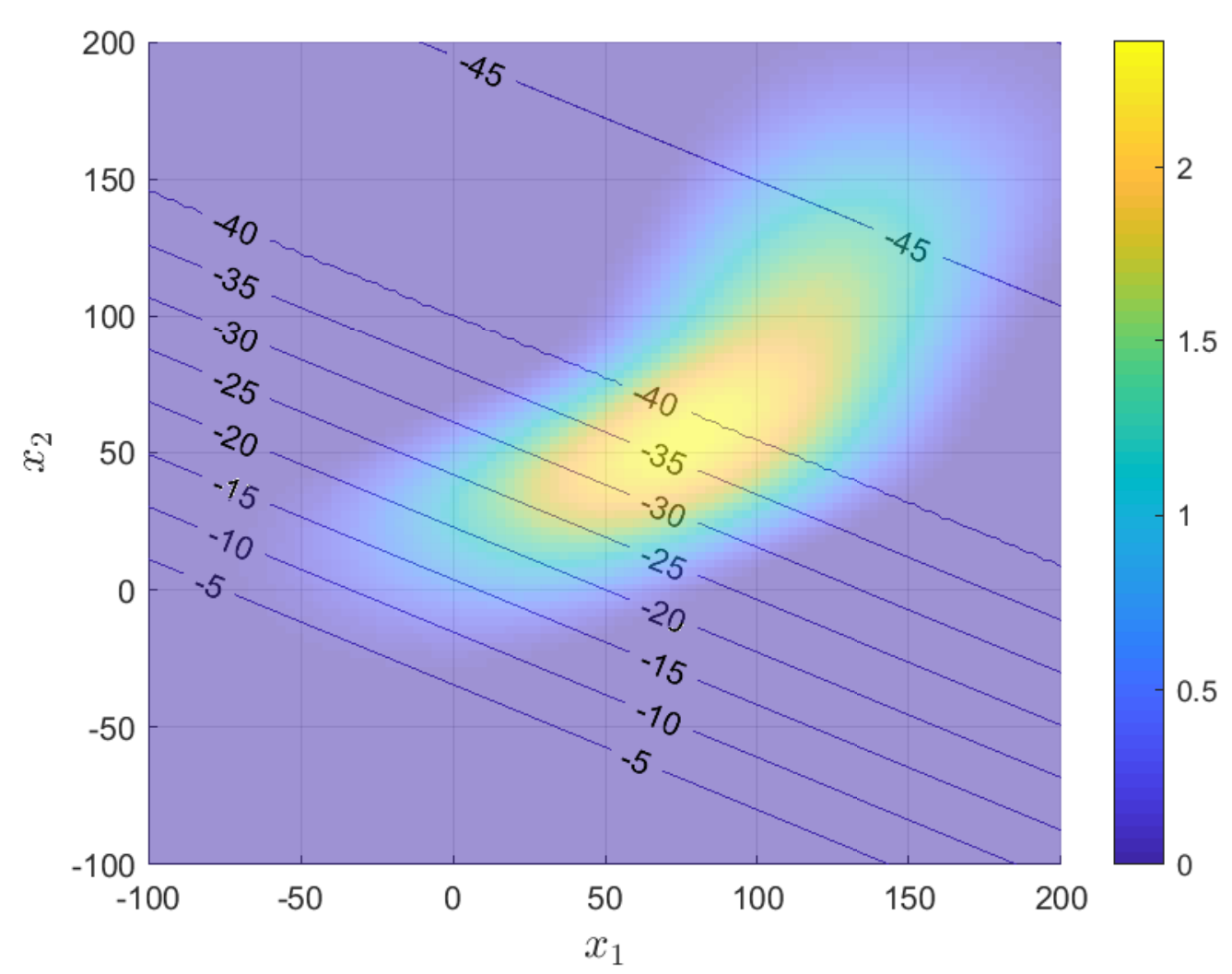} 
    \caption{The slip field $\cal{G}$. We are using a pure thrust model, that is,
		 the direction of the slip is in the line of steepest descent,  only the modulus of 
		$\cal{G}$ is shown. The colorbar shows the scale for the modulus of $\cal{G}$.
		As previously, lines of equal depth on $\Gamma$ are indicated.}
    \label{the_slip}
\end{figure}

\begin{figure}[htbp]
   \centering
        \includegraphics[scale=.5]{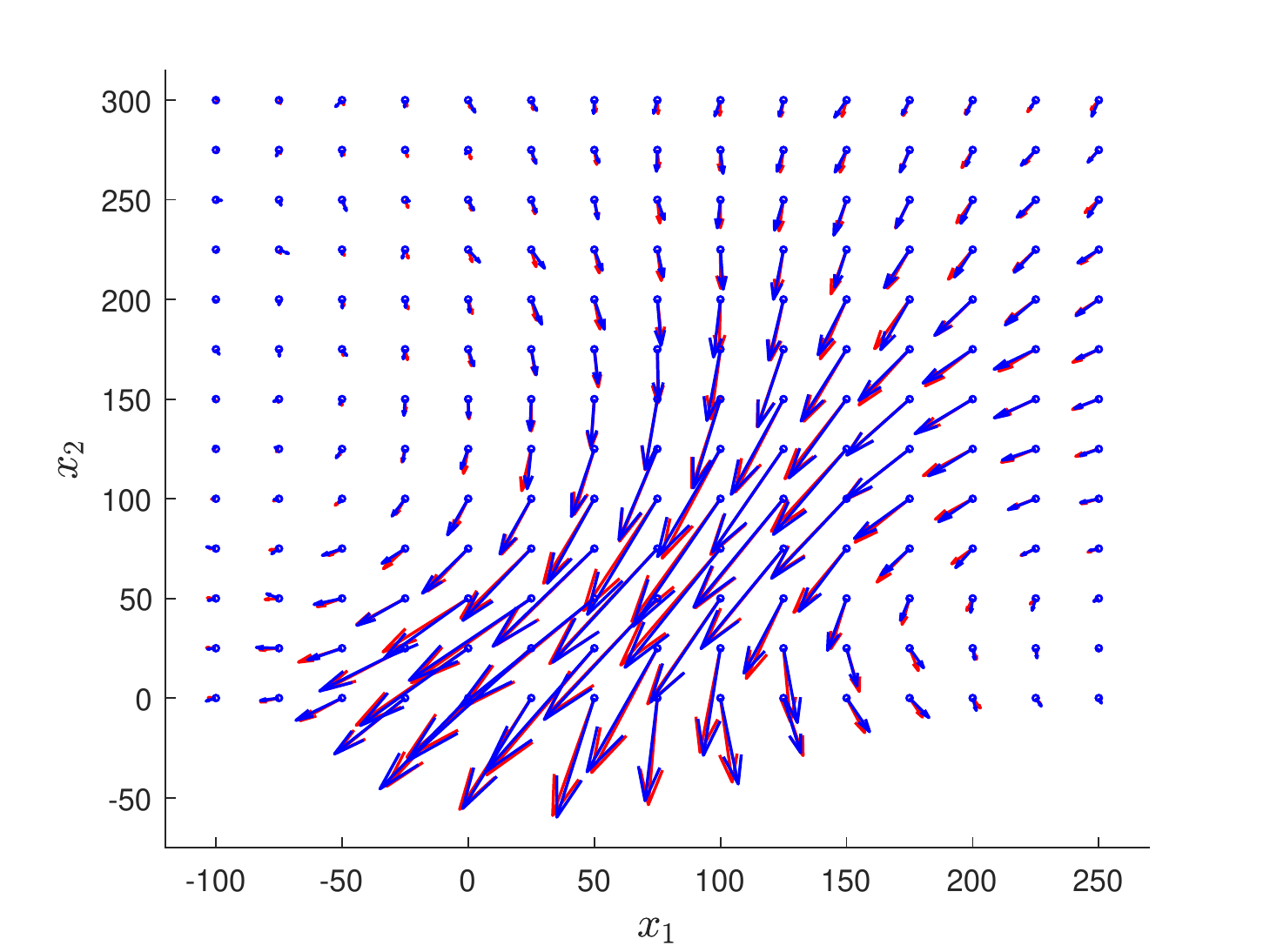} 
        \includegraphics[scale=.5]{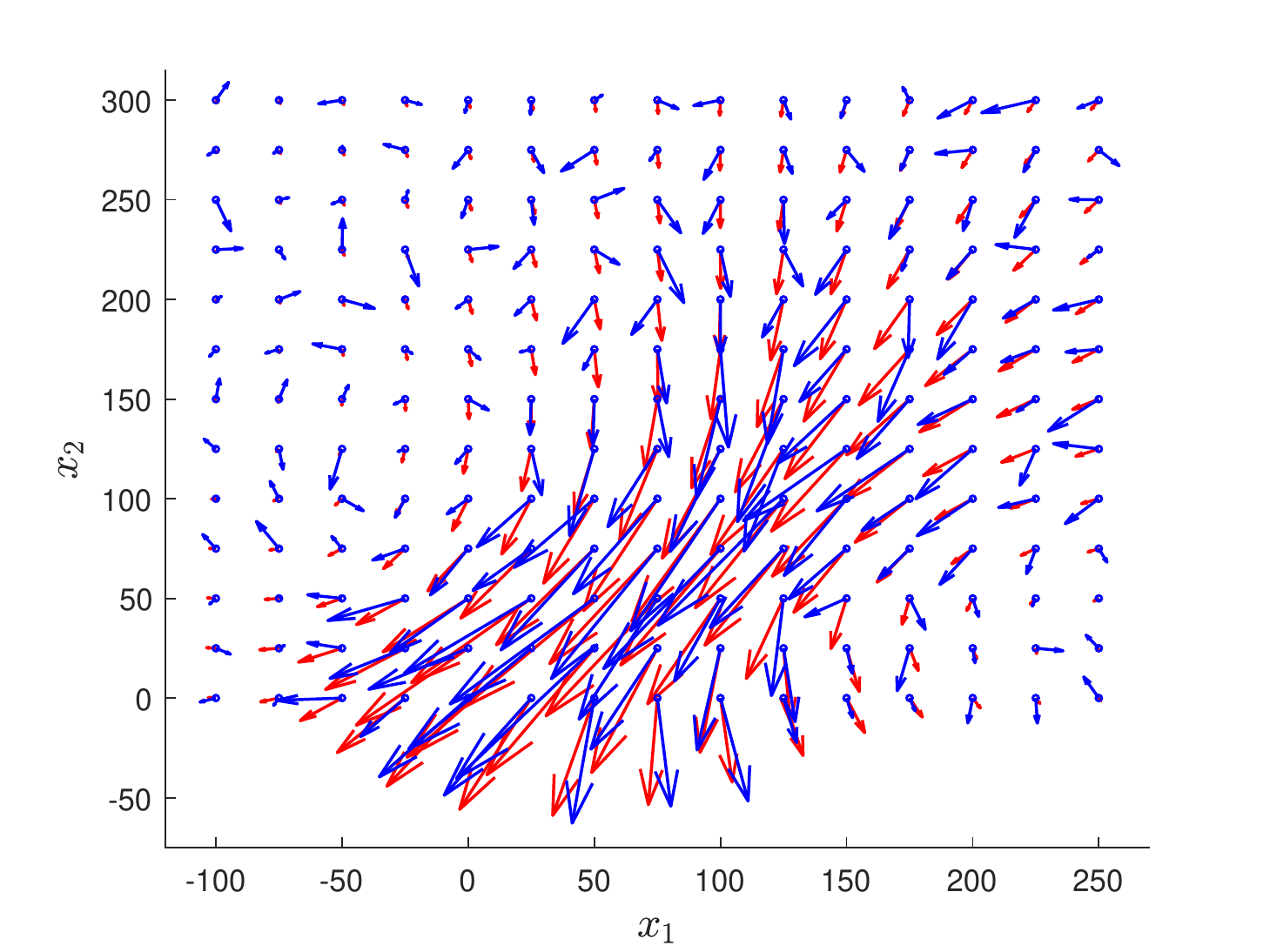} 
    \caption{Plots of the data for solving the fault inverse problem in  our simulations.
		The data  are
		realizations of noisy horizontal surface displacements (blue) at the measurement points obtained
		by using the slip field shown in Figure \ref{the_slip} occurring on the fault $\Gamma$ 
		shown in Figure \ref{fault}. 
		Top: low noise scenario. Bottom: high noise scenario. 
		In each graph the unperturbed field is sketched in red.}
    \label{surface_disp}
\end{figure}

We consider data generated in a configuration
 closely related to  studies involving field data for a particular region and a specific seismic event 
\cite{volkov2019stochastic, volkov2017determining}.
In  those studies, 
simulations 
involved
only planar
faults,  while here we examine  the case of
fault geometries defined by pairs of contiguous quadrilaterals.
Obviously, reconstructing finer geometries as considered here requires many
more measurement points than used in  \cite{volkov2019stochastic, volkov2017determining}
(11 points in \cite{volkov2019stochastic} versus 195 points here).
The higher number of measurement points used here allows
us to reconstruct $\Gamma$
even if the data is very noisy, at the cost of
finding large standard deviations.
In our model, the geometry of $\Gamma$ is determined from $\bm$ 
in $\RR^6$ in the following way:
\begin{itemize}
\item $\Gamma$ is beneath the square $[-100,200] \times [-100, 200]$
\item Let $P_1$ be the point  $(-100,-100, m_1)$
\item Let $P_2 $ be the point $(-100, m_2, m_3)$, such that $ -100 <m_2 <  200$
\item Let $P_3 $ be the point $(200, m_4, m_5)$, such that $ -100 <m_4 <  200$
\item Let $P_4 $ be the point $(200, 200, m_6)$
\item Let $P_5$ be the point in the plane $P_1 P_2 P_3$ with $x_1, x_2$ 
coordinates $(200, -100)$
\item Let $P_6$ be the point in the plane $P_2 P_3 P_4$ with $x_1, x_2$ 
coordinates $(-100, 200)$ 
\item Form the union of the two quadrilaterals $P_1 P_2 P_3 P_5$ 
and $P_2 P_3 P_4 P_6$ and discard the part where $x_3 \geq 0$
to obtain $\Gamma$
\end{itemize}
For generating forward data we picked the particular values
\bean
 m =
(24,     145,    -40,     8,    -40   -50). \label{values}
\eean
 In Figure \ref{fault}, we show a sketch of  $\Gamma$ viewed from above with the points 
$P_1, P_2, P_3, P_4, P_5, P_6$, and  we also sketch $\Gamma$ 
in three dimensions.
In Figure \ref{the_slip}, we show a graph of  the slip field 
${\cal G}$ as a function of $(x_1, x_2)$ (recall that this slip field is supported
on $\Gamma$ so Figure \ref{the_slip} shows a projection of ${\cal G}$  on a horizontal plane).
We model a slip of pure thrust type, meaning that slip occurs in the direction of steepest descent,
so only the norm of ${\cal G}$   is graphed.
We used this slip $\cal{G}$  to compute the resulting surface displacements 
thanks to formula (\ref{int formula})  
where we discretized the integral on a fine mesh.
The data for the inverse problem is the noise-free three dimensional displacements
$\bu_{free}$
 at the measurement points shown in 
Figure \ref{fault}  to which we added
 Gaussian noise with zero mean and covariance $\sigma^2 I$ to obtain $\bu$.
We consider two scenarios: lower and higher noise.
In the lower noise  scenario $\sigma$ was set to be equal to 5\%
of the maximum of the absolute values of the components of $\bu_{free}$
(in other words, 5\% of $\| \bu_{free} \|_{\infty}$).
For the particular realization used in solving the inverse problem, 
this led to a relative error in Euclidean norm of about 7\%
(in other words, for this particular realization ${\cal E} (\omega) $ of the noise,
$\| {\cal E} (\omega) \|/\| \bu_{free} \|$ was about 0.07, where $\| . \| $ is the Euclidean norm).
In the higher noise case scenario $\sigma$ was set to be equal to 25\%
of the maximum of the absolute values of the components of $\bu$
(in other words, 25\% of $\| \bu_{free} \|_{\infty}$).
This time, this led to 
 a relative error in Euclidean norm of about 37\%.
Both realizations are shown in Figure \ref{surface_disp}
(only the horizontal components are sketched for the sake of brevity).
All lengthscales used in these simulation 
are in line
with  
the canonical example  from  geophysics provided by the 
2007 Guerrero slow slip event \cite{volkov2019stochastic, volkov2017determining}.
In particular,
$x_1, x_2, x_3$ are thought of as  given in kilometers and  $\bg$ and $\bu$  in meters.
In real life applications, the noise level is likely to lie somewhere between the low noise scenario
and high noise scenario considered here. 

\subsection{Numerical results from our  parallel sampling algorithm} \label{our results}
\begin{figure}[htbp]
   \centering
        \includegraphics[scale=.4]{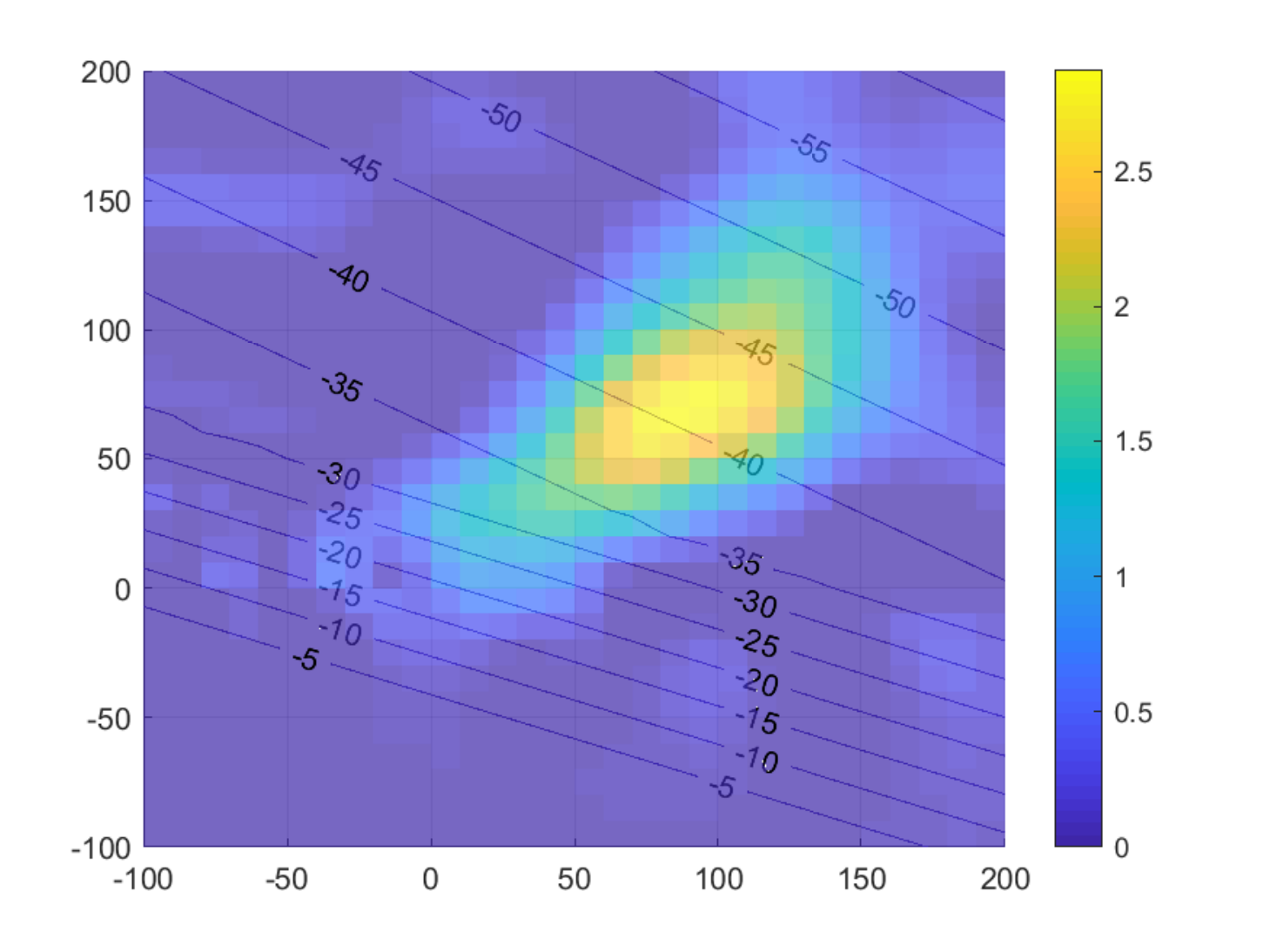} 
        \includegraphics[scale=.4]{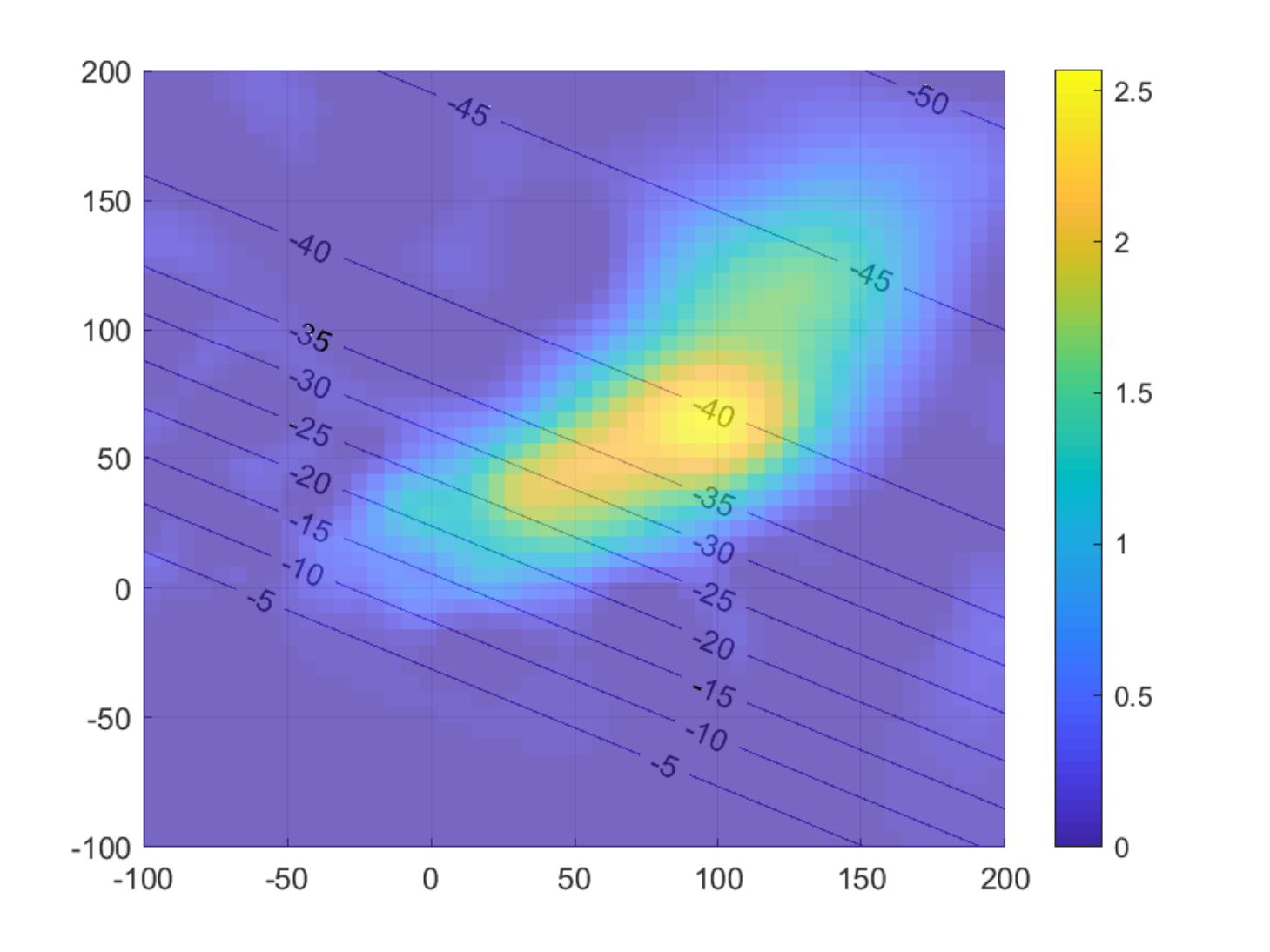} 
        \includegraphics[scale=.4]{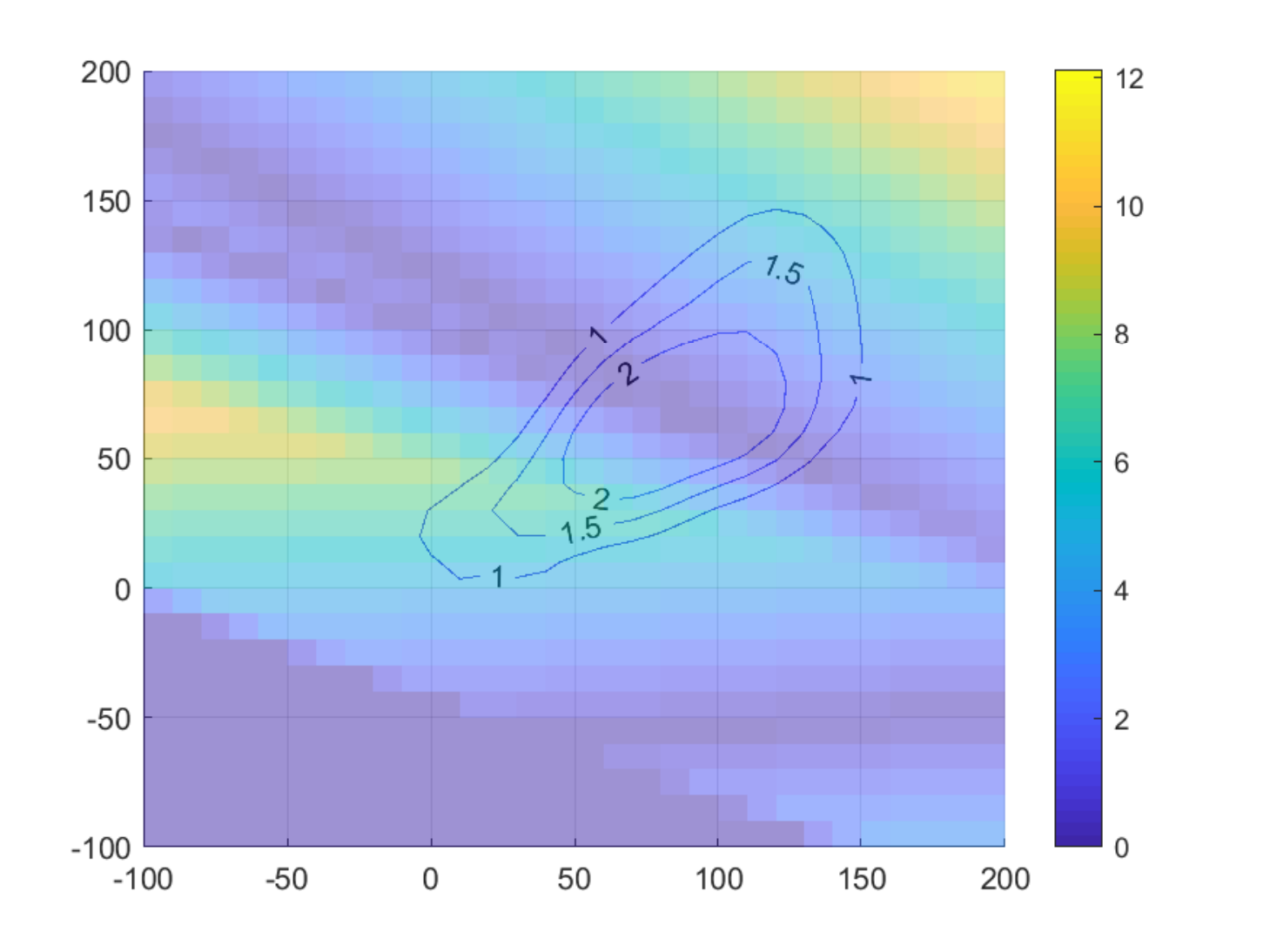} 
        \includegraphics[scale=.4]{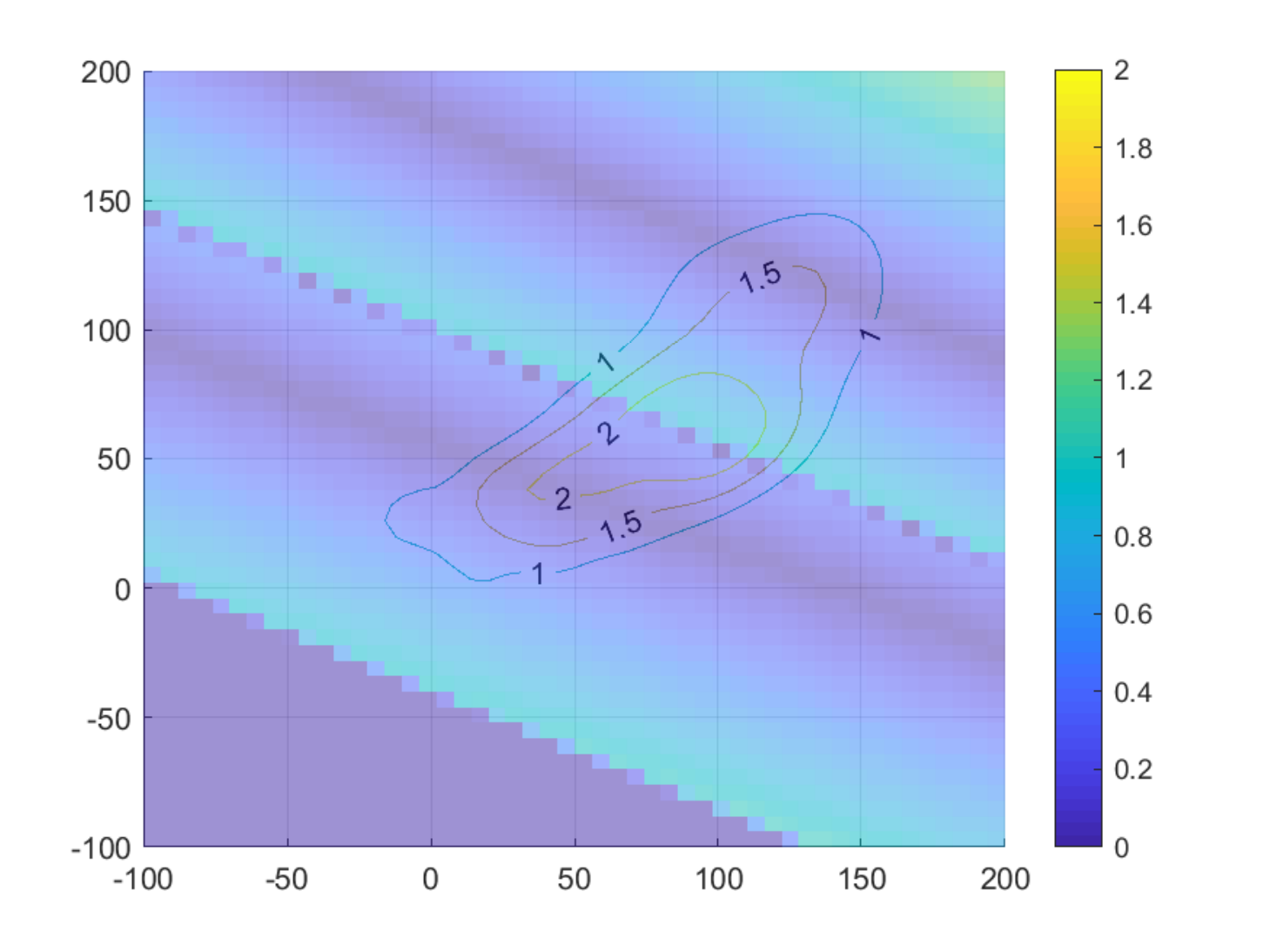} 
        \includegraphics[scale=.4]{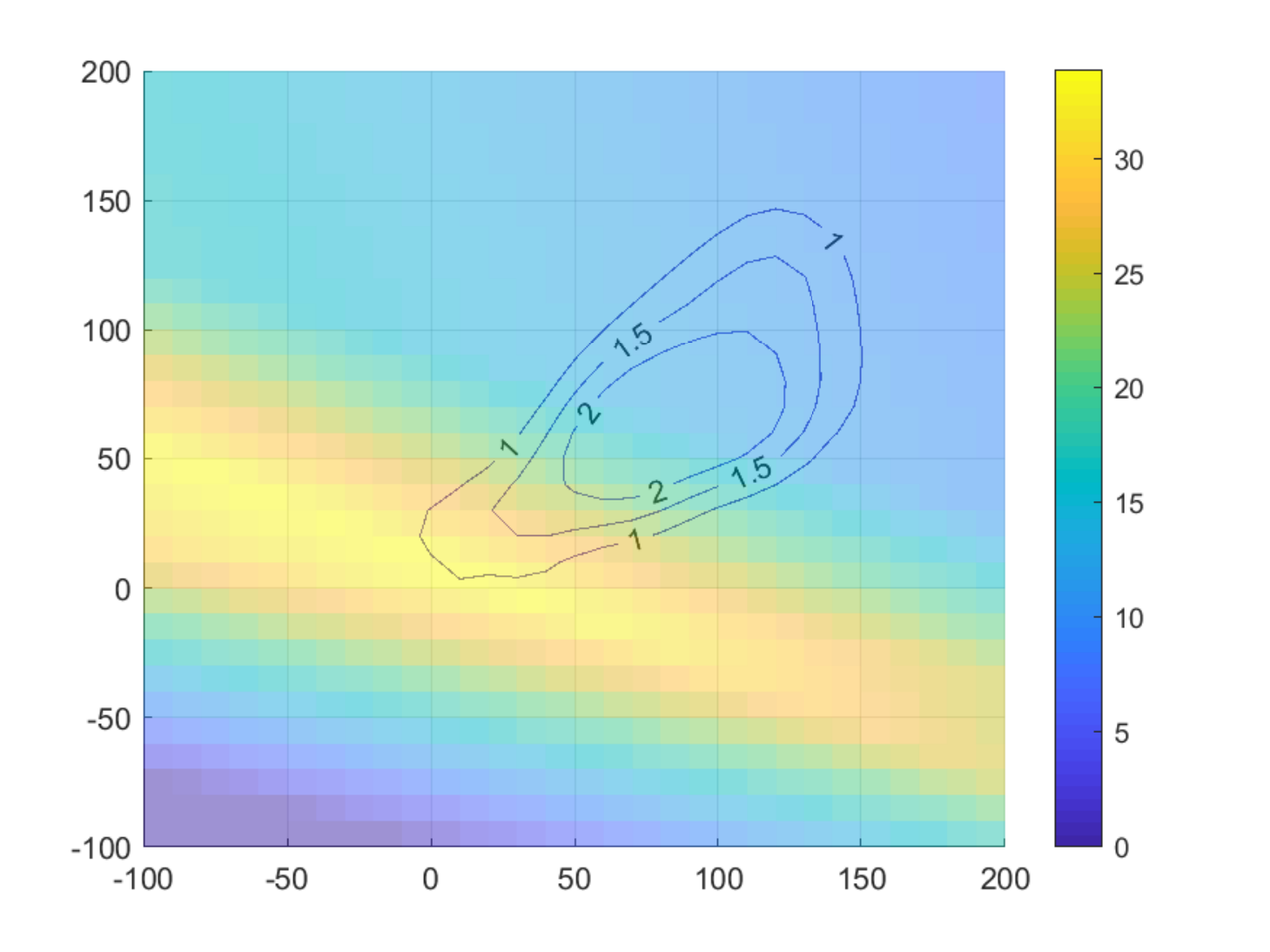} 
        \includegraphics[scale=.4]{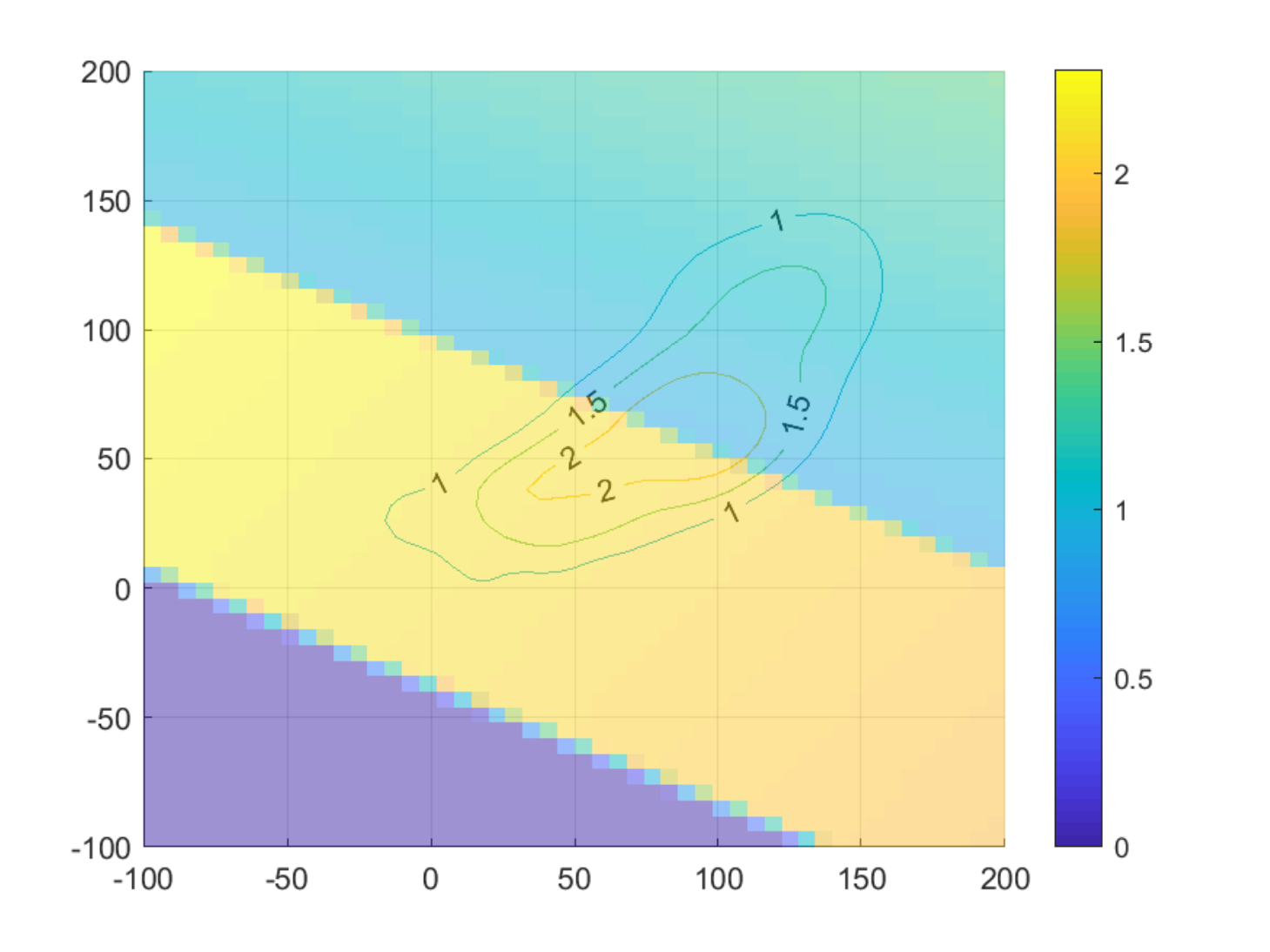} 
    \caption{
		Reconstructed slip field ${\cal G}$ (discretized by $\bg_{min}$) and computed
		expected depth profile for the fault
		$\Gamma$ 
		reconstructed
		from the data shown in Figure \ref{surface_disp}. 
		In all six graphs in this figure, the horizontal axis is for $x_1$, and the vertical axis
		is for $x_2$.
		Left column: high noise scenario. Right column: 
		low noise scenario. 
		First row: depth contour lines for $\Gamma$ corresponding to the expected value of 
		$\bm$ and the slip field modulus (shown in color) given that geometry and the computed expected value of $\alpha$. 
		The slip field was  obtained
		by applying  the formula 
		$\bg_{min} = (A_\bm' A_\bm + \alpha R'R)^{-1} A_\bm' \bu$.
		Second row: absolute value of the difference between computed expected
		depth of $\Gamma$ and true depth as a function 
		of $(x_1,x_2)$
		shown in color,   with contour lines of reconstructed 
	$|{\cal G}|$.
	The third  row shows 
	two standard deviations for the reconstructed depth as a function of $(x_1, x_2)$, 
	 with again contour lines of reconstructed 
	$|{\cal G}|$.
		}
    \label{errors_stds}
\end{figure}

\begin{figure}[htbp]
   \centering
        \includegraphics[scale=.4]{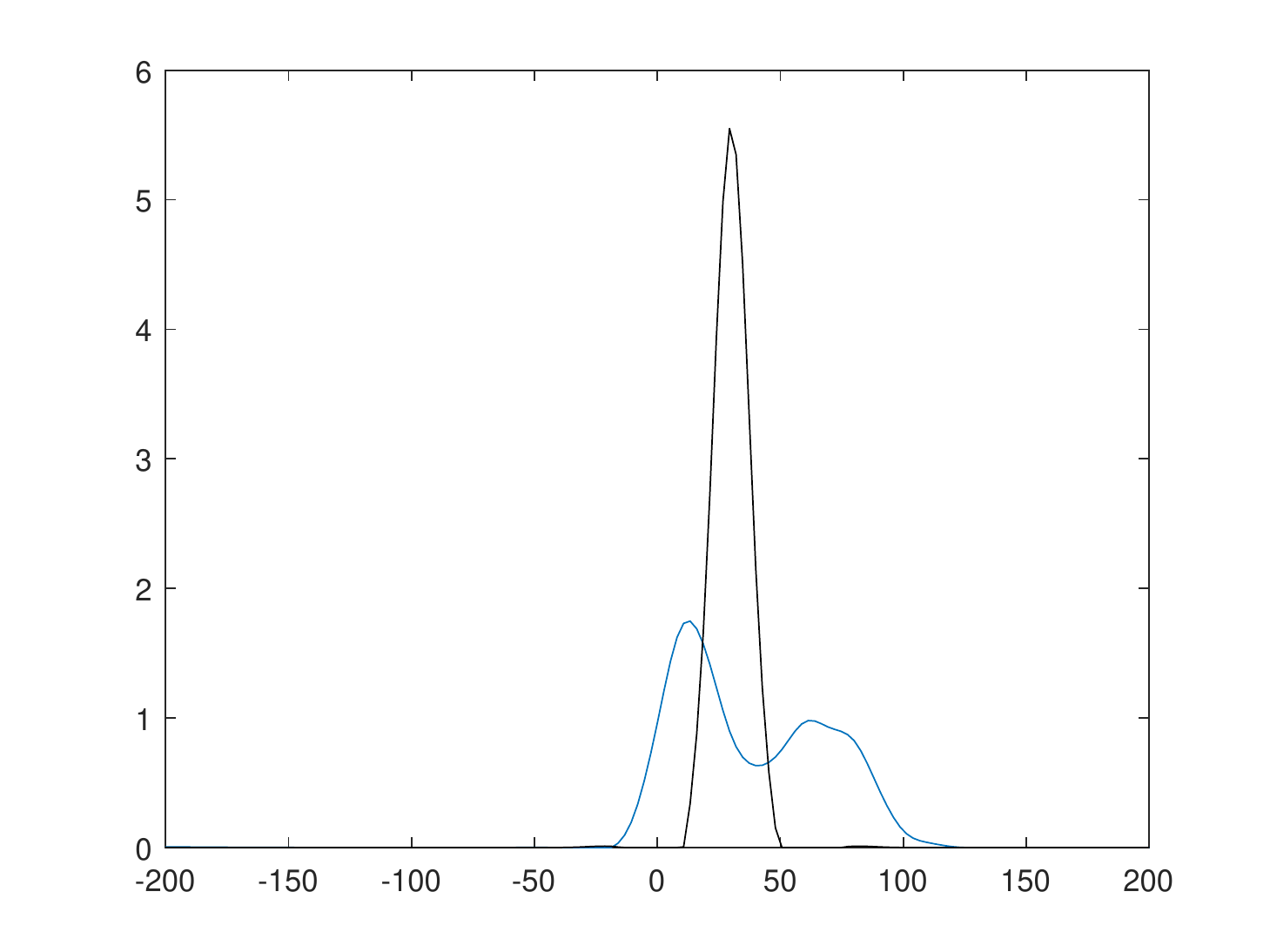} 
        \includegraphics[scale=.4]{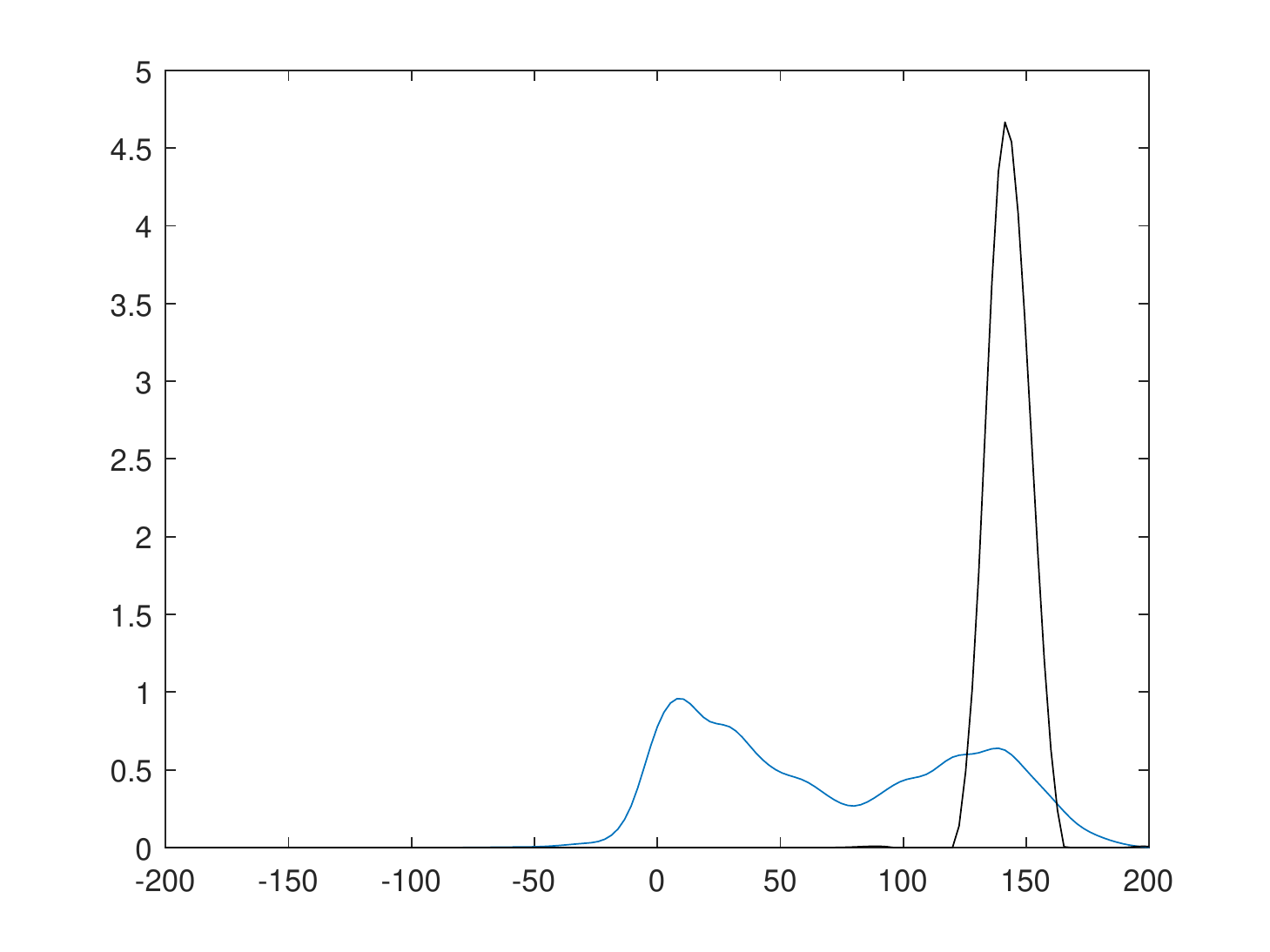} 
        \includegraphics[scale=.4]{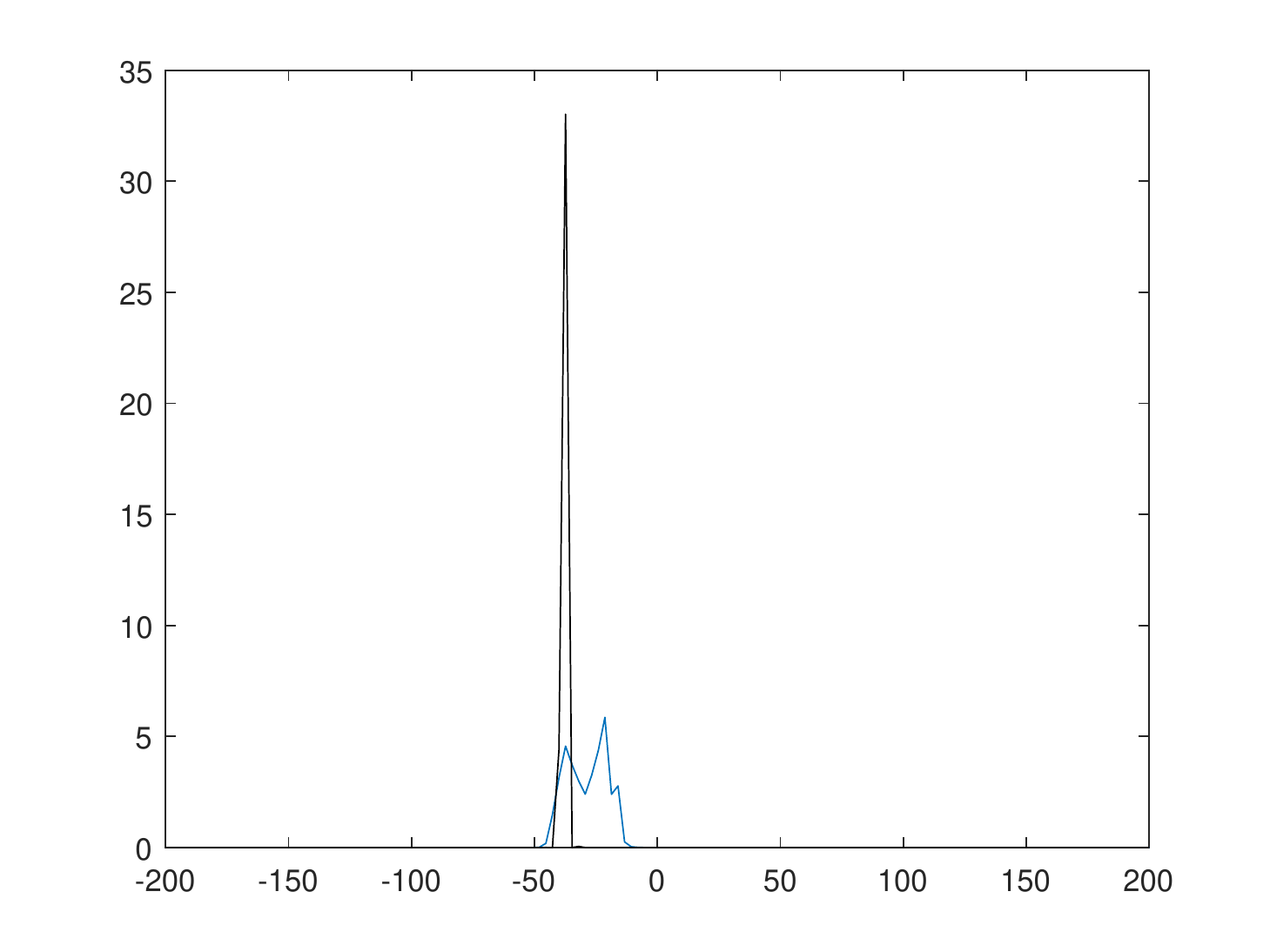} 
        \includegraphics[scale=.4]{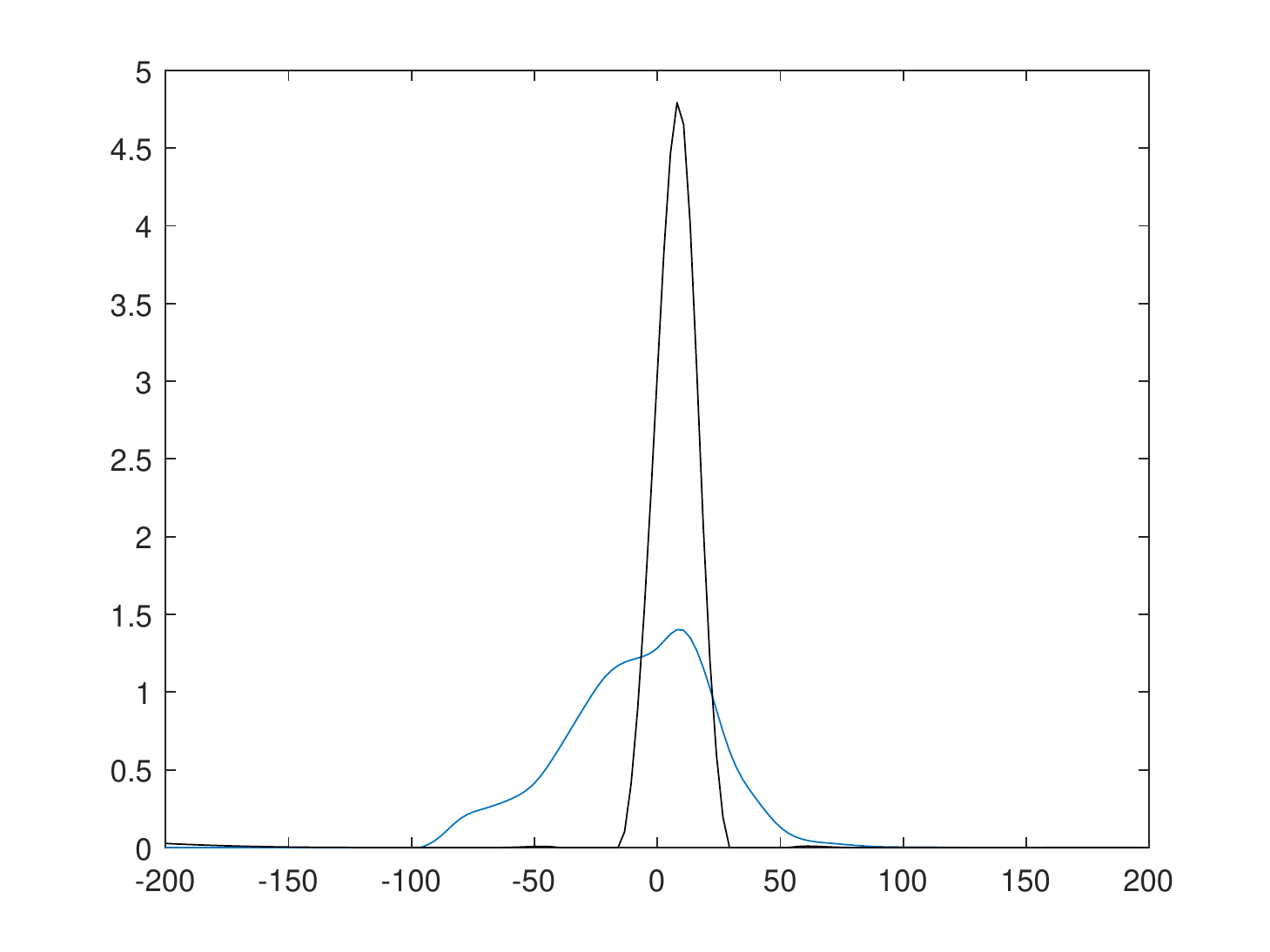} 
        \includegraphics[scale=.4]{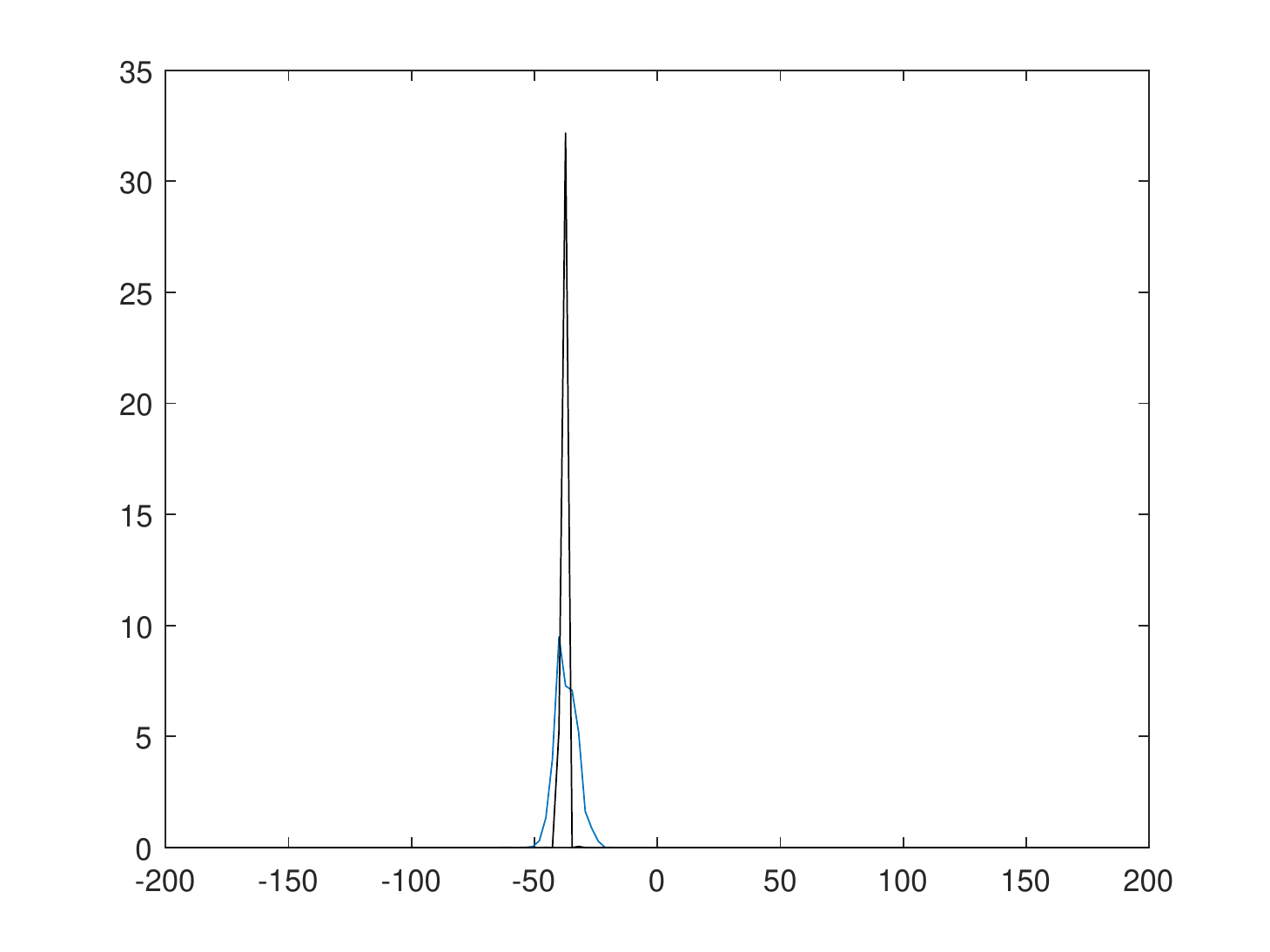} 
        \includegraphics[scale=.4]{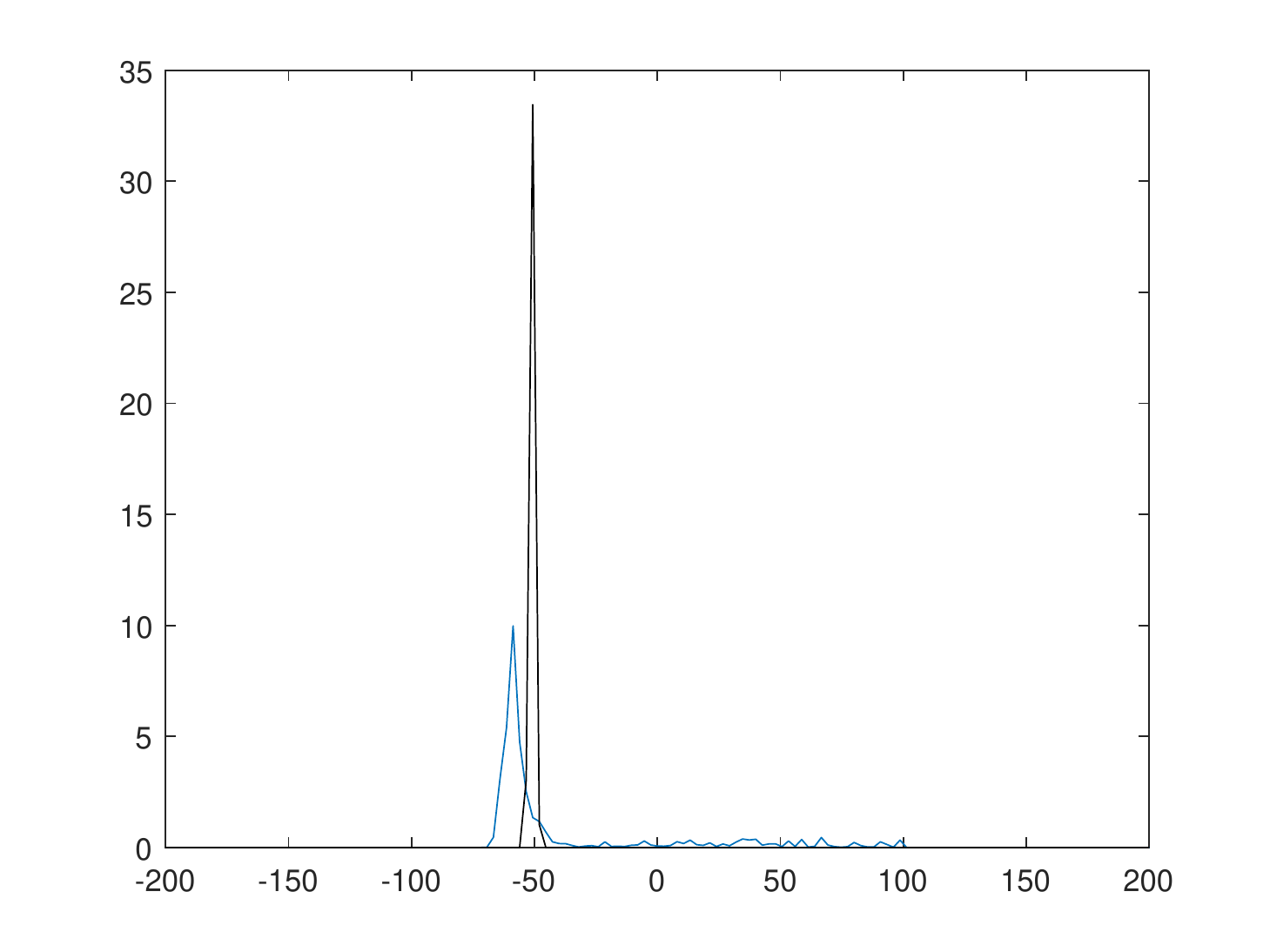} 
			 \includegraphics[scale=.4]{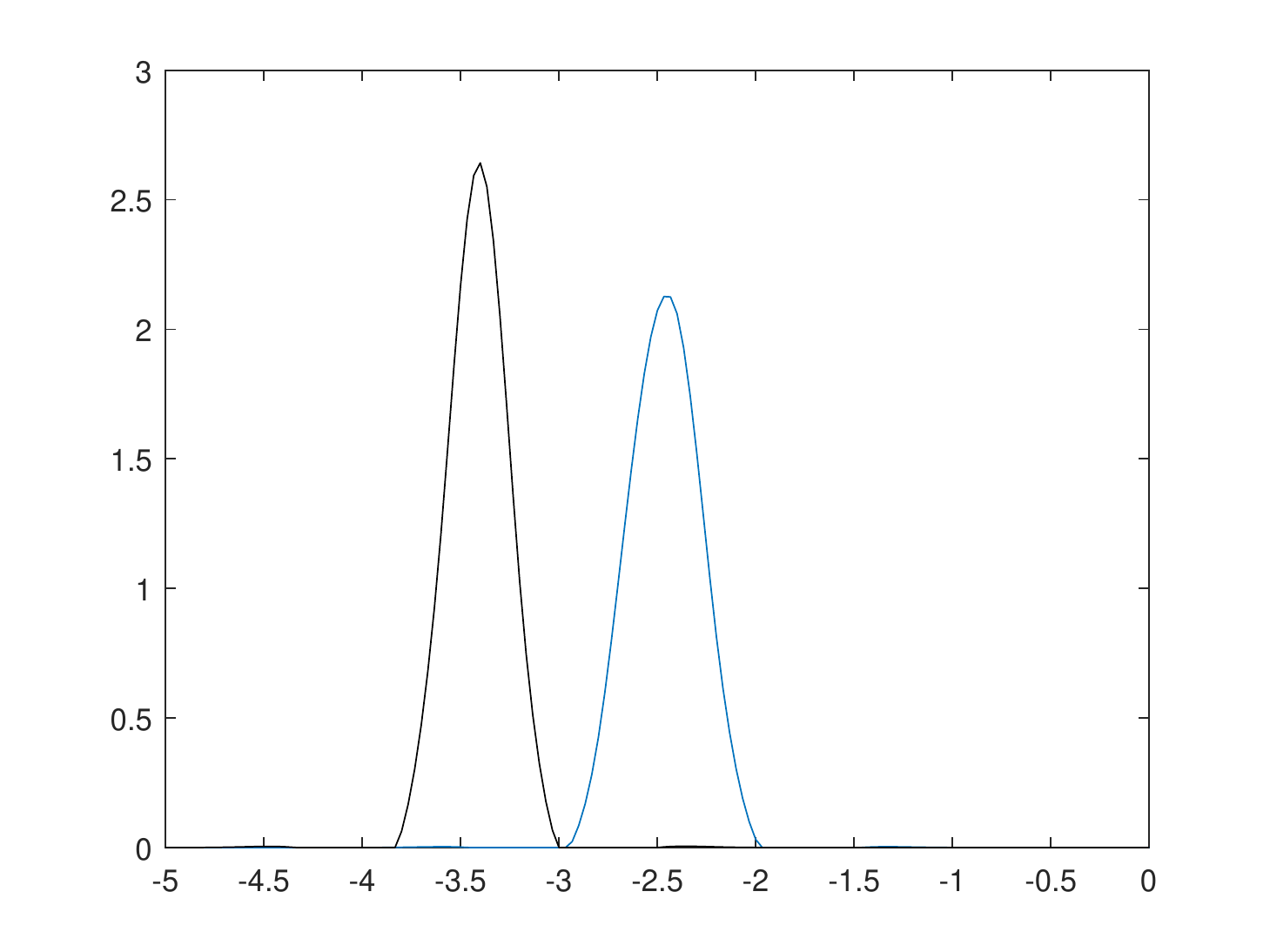} 
    \caption{
		Reconstructed posterior marginal distribution functions for the six components of $\bm$
		and for $\log_{10} \alpha$. Blue: high noise scenario. Black: low noise scenario.
		For $\bm$, true  values are nearly indistinguishable from the peak of the black curves.}
    \label{pdfs}
\end{figure}

\begin{figure}[htbp]
   \centering
        \includegraphics[scale=.4]{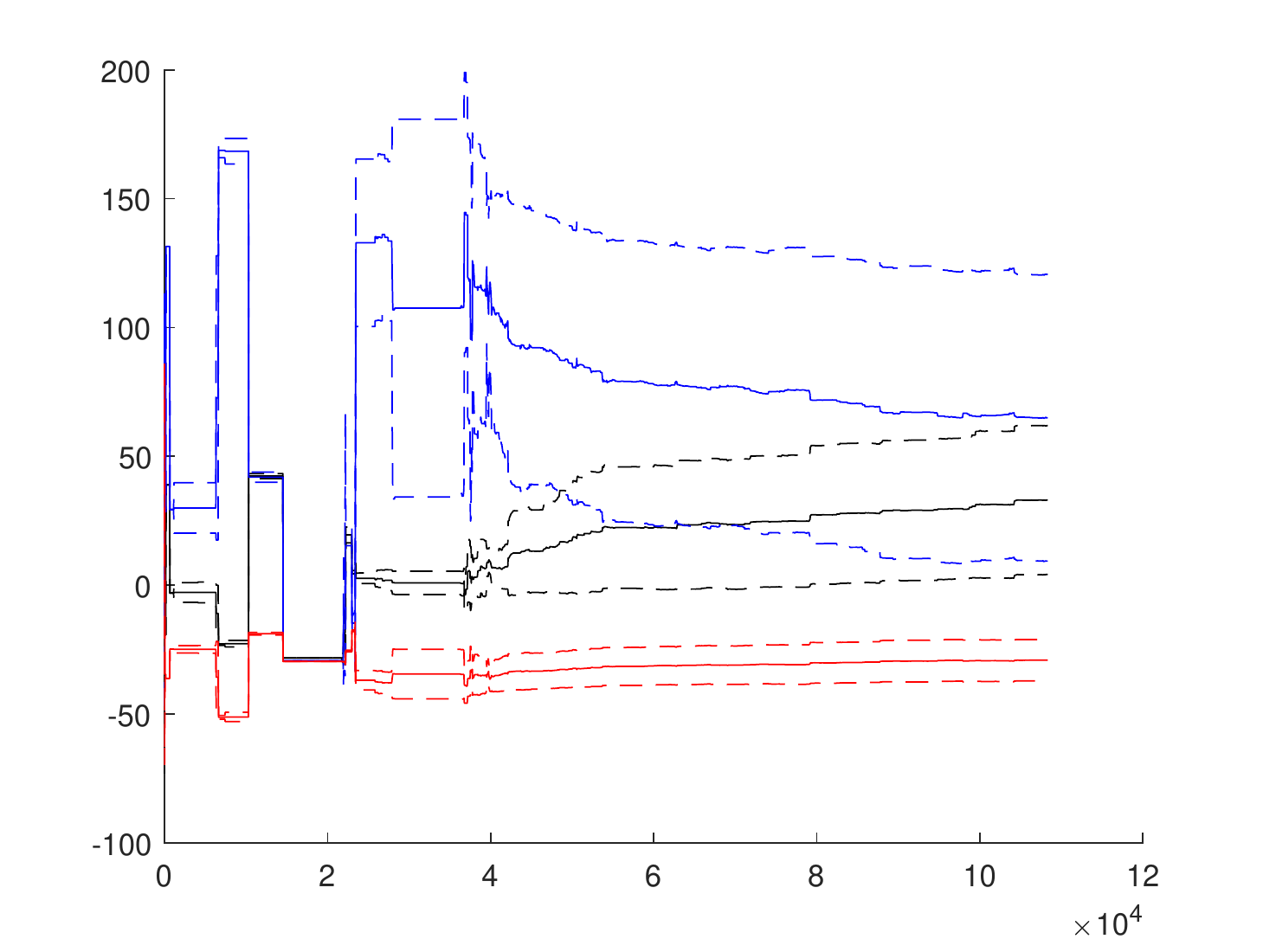} 
				        \includegraphics[scale=.4]{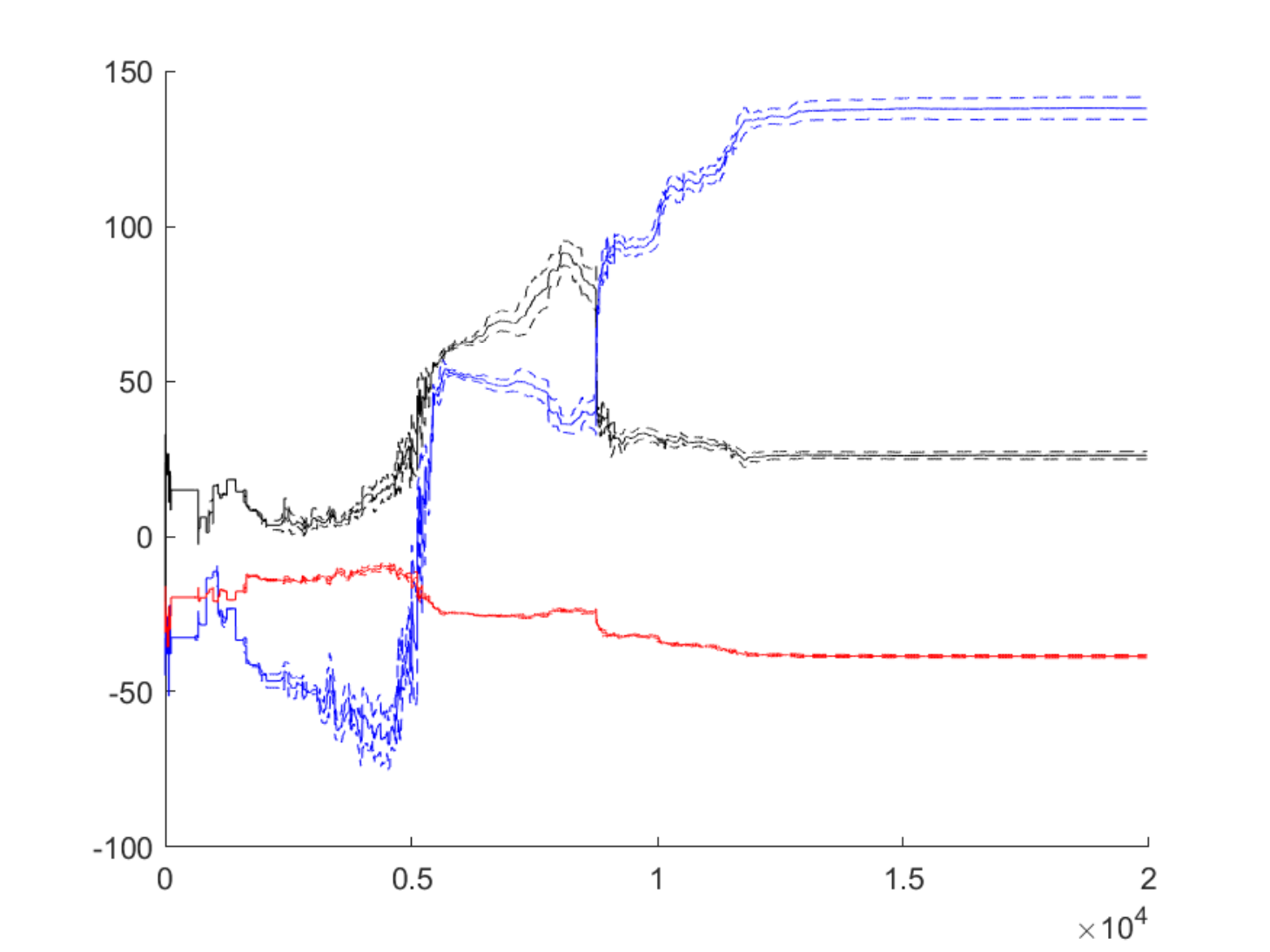} 
        \includegraphics[scale=.4]{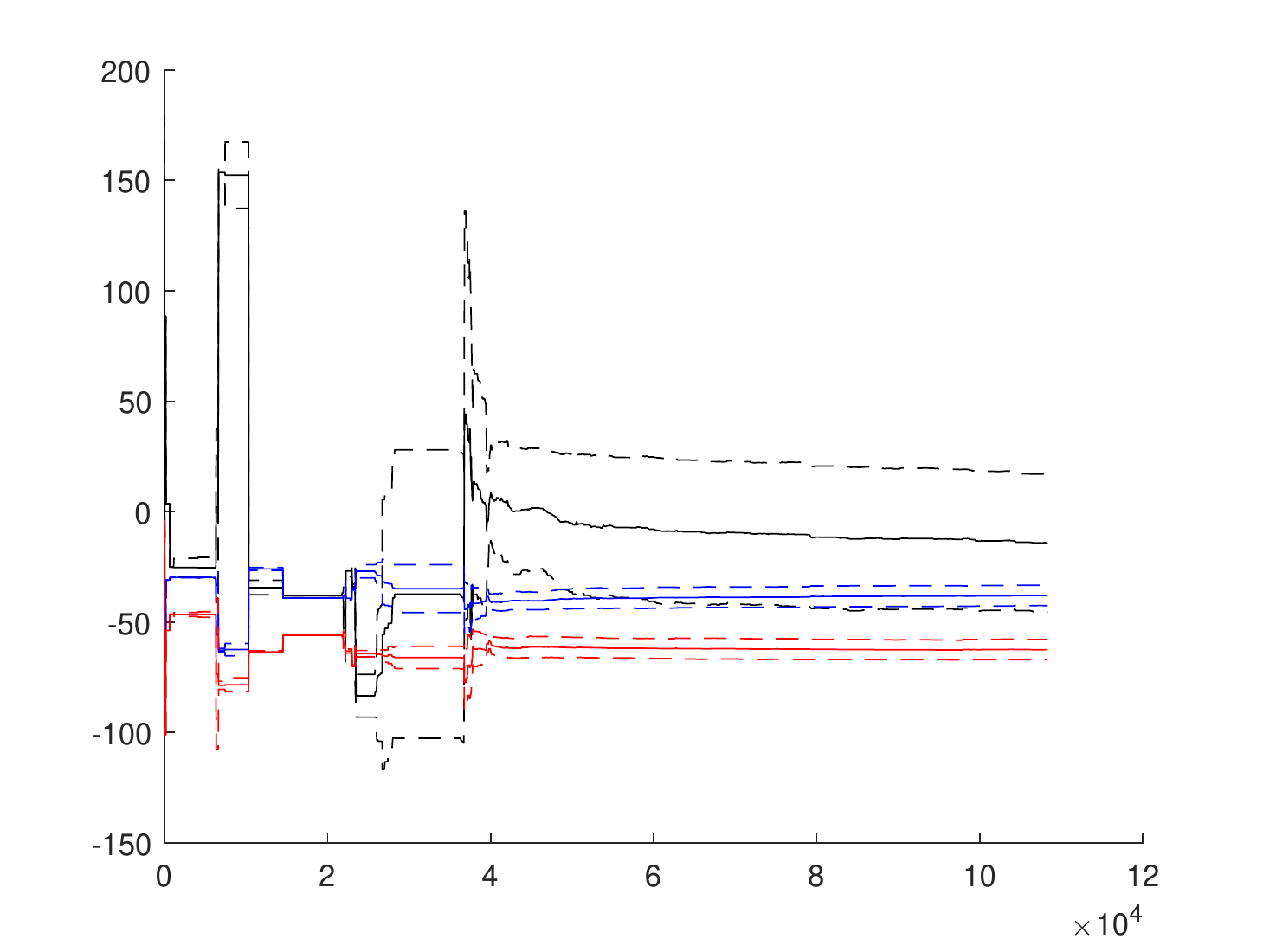} 
				        \includegraphics[scale=.4]{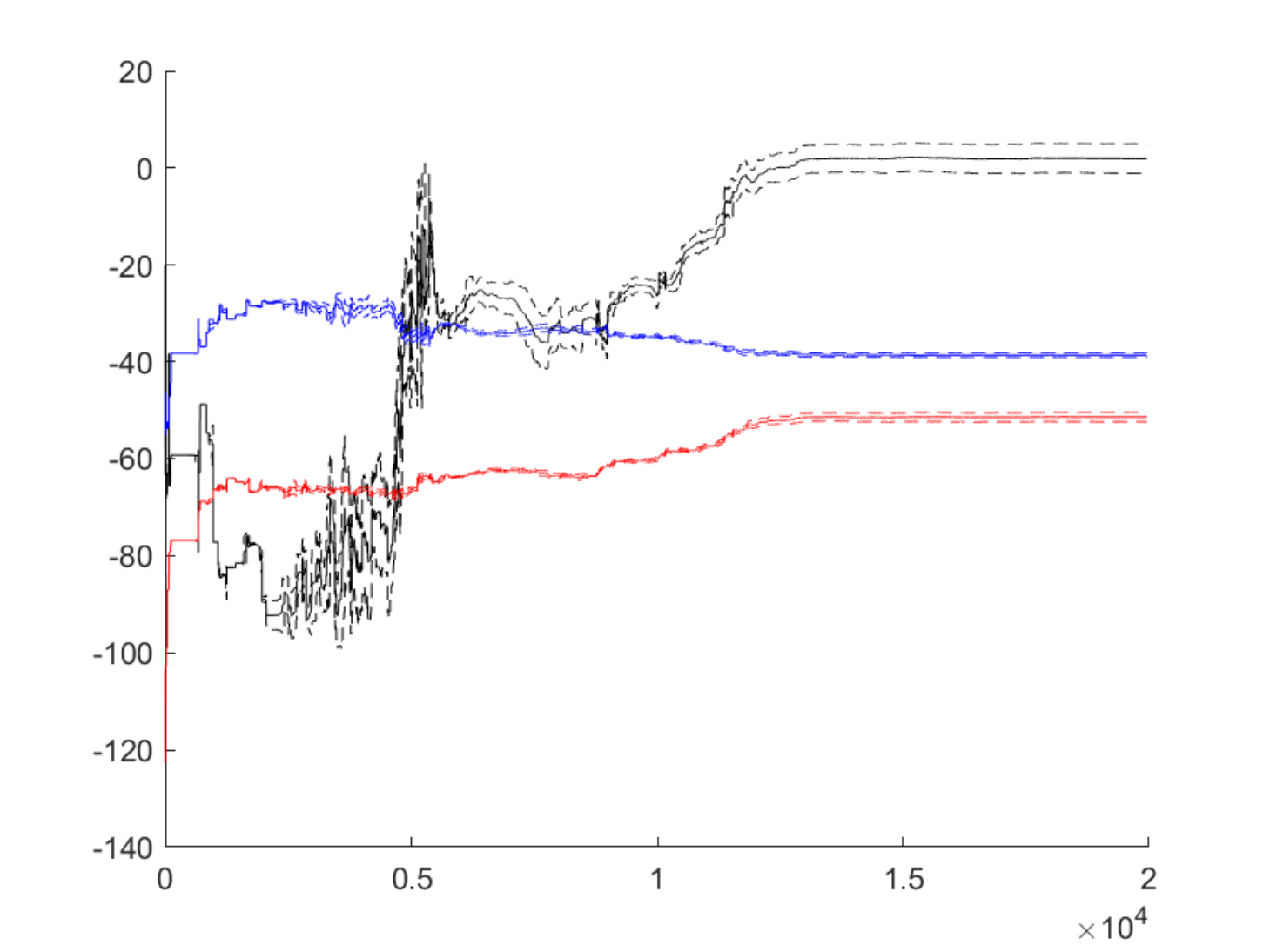} 
    \caption{
	First row: evolution of computed expected values of $m_1,m_2,m_3$ (black, blue, and red solid lines) and
	one standard deviation envelopes. Left: high noise scenario. Right: low noise scenario.
	Second row: evolution of computed expected values of $m_4,m_5,m_6$ (black, blue, and red solid lines) and
	one standard deviation envelopes. Left: high noise scenario. Right: low noise scenario.}
    \label{progress}
\end{figure}

The parallel sampling algorithm introduced in section 
\ref{par} requires the  knowledge of  a prior distribution for the random variable $(\bm, \alpha)$.
Here,
we assume that the priors of $\bm$ and $\alpha$ are independent.
The prior of $\bm$ is  chosen to
be uniformly distributed on the  subset ${\cal B}$ of  $[-200,200]^6$ 
such that the angle $\theta$ between  the    vectors normal  to the two
quadrilaterals whose union is $\Gamma$ 
satisfies  $\cos \theta  \geq 0.8$. That way, the angle between these two 
quadrilaterals is between 143 and 217 degrees: this can be interpreted as a regularity condition
on the slip field $\bg$ since it was set to point in the direction of steepest descent.
As to $\alpha$, we assumed that $\log_{10} \alpha$ follows a uniform prior on $[-5, 0]$. \\
Next, we  present  results obtained by applying 
our parallel sampling  algorithm 
to the data shown in Figure \ref{surface_disp} for  the low noise and the high noise scenario.
Computations were performed on a parallel platform that uses $N_{par}=20$
processors.
After computing the expected value of $\bm$, we sketched the 
corresponding geometry of $\Gamma$
(Figure \ref{errors_stds}, first row) by plotting depth contour lines.
On the same graph we plotted the magnitude of the reconstructed  slip field on $\Gamma$ as a function
of $(x_1, x_2)$.  This reconstruction was done by using the expected value of 
$(\bm,\alpha)$ and solving the linear system $(A_\bm' A_\bm + \alpha R'R)\bg_{min} =  A_\bm' \bu$.
The reconstructed slip field 
is very close to the true one shown in Figure  \ref{the_slip}.
In the second row of  Figure \ref{errors_stds}, we show the absolute value of the
difference between the computed  expected depth $x_3$ on   $\Gamma$
minus its true value
%
and a contour profile of the reconstructed slip $\bg_{min}$. 
It is  noticeable  that the error in reconstructing the depth in the high noise scenario
is rather lower where $\bg_{min}$ is larger: this is in line
with previous theoretical studies \cite{volkov2019stochastic, volkov2017reconstruction}
where it was proved that fault geometries can only  be reconstructed on the support of
slip fields. 
In the third row of Figure \ref{errors_stds},
we show
two standard deviations for the reconstructed depth as a function of $(x_1, x_2)$, 
	 with again contour lines of 
	$\bg_{min}$.
We note that this difference (close to 2) is very low in the low noise scenario compared to the depth
at the center of the support of $\bg_{min}$ (close to -40).
In the high noise scenario
the standard deviation is rather lower in the region where $|\bg_{min}|$ is larger.
We show in  Figure \ref{pdfs} 
reconstructed posterior marginal distribution functions for the six components of $\bm$
		and for $\alpha$.  
		Interestingly, we notice that the range of high probability for $\alpha$ is much higher
		in the higher noise scenario
		(more than 10 times higher, since the graph is that of $\log_{10} \alpha$).
		Intuitively, it is clear that stronger noise would require 
		more regularization, as  Morozov principle dictates \cite{vogel2002computational},
		but the strength of our algorithm is that it automatically selects a 
		good range for $\alpha$ without user input or prior knowledge about  $\sigma$. 
		In the low noise scenario, for all six components of $\bm$, the
		reconstructed posterior marginal distribution functions 
		peak very close to their true value, the difference would not actually be visible on the 
		graphs. 
		The picture is quite different in the high noise scenario. 
		The support of the distribution functions are much wider in that case
		and two peaks are apparent for $m_1$ and for $m_2$.
		The large width for these distribution functions  is related
		to a much larger number of samples in the random walk for the algorithm to converge 
		as illustrated in Figure \ref{progress}.
		As to choosing a starting 
 point $(\bm_1, \alpha_1)$, 
		we found that  it  was most efficient to  draw $N_{burn} \times N_{par}$ 
		 samples from the prior and  use these samples to compute an expected value,
		which we set to be equal to $(\bm_1, \alpha_1)$.\\
		To conclude this section, we would like to emphasize that the numerical results that 
		we show in this paper are not so sensitive to the the particular realizations of 
		the noise and the intrinsic randomness of Markov chains. 
			We conducted a large number of simulations each starting from
			a different realization of $\bu$. In the low noise scenario, the differences between
			 estimates of expected values and covariances were negligible.
			In the high noise scenario, the differences between
			 estimates of expected values and covariances were more  appreciable,
			however  the difference between the expected value of the depth of $\Gamma$
			and its true value
			as in row 2 of Figure \ref{errors_stds}, and the plot for  the  
			standard deviation  of reconstructed depth as in row 3 of Figure \ref{errors_stds},
			were comparable and the differences were small in the region of high values
			of the reconstructed slip $\bg_{min}$.

%

\subsection{Comparison to methods based on  GCV, ML, and CLS}\label{comparisons}
\begin{figure}
   \centering
        \includegraphics[scale=.4]{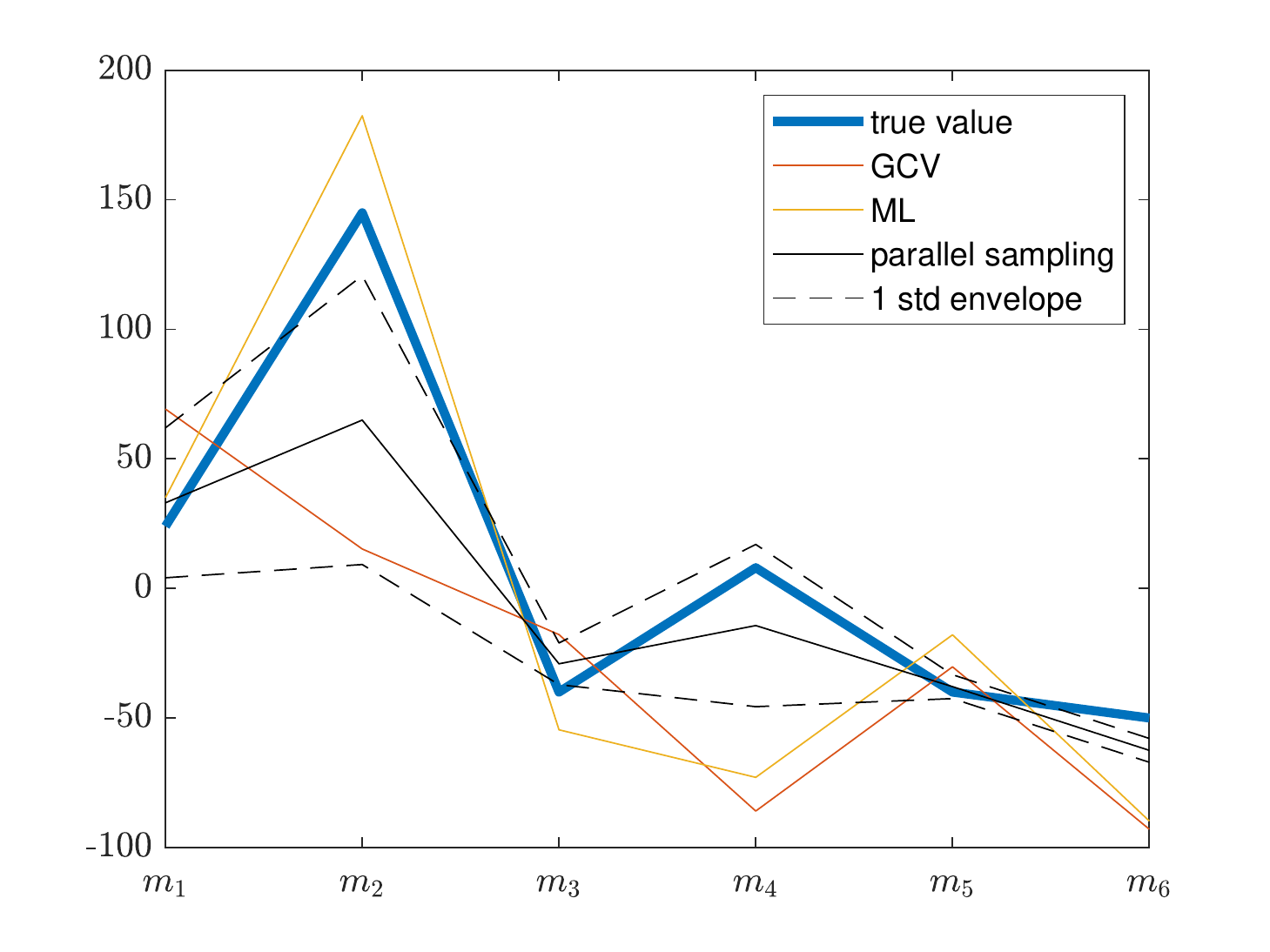} 
			  \includegraphics[scale=.4]{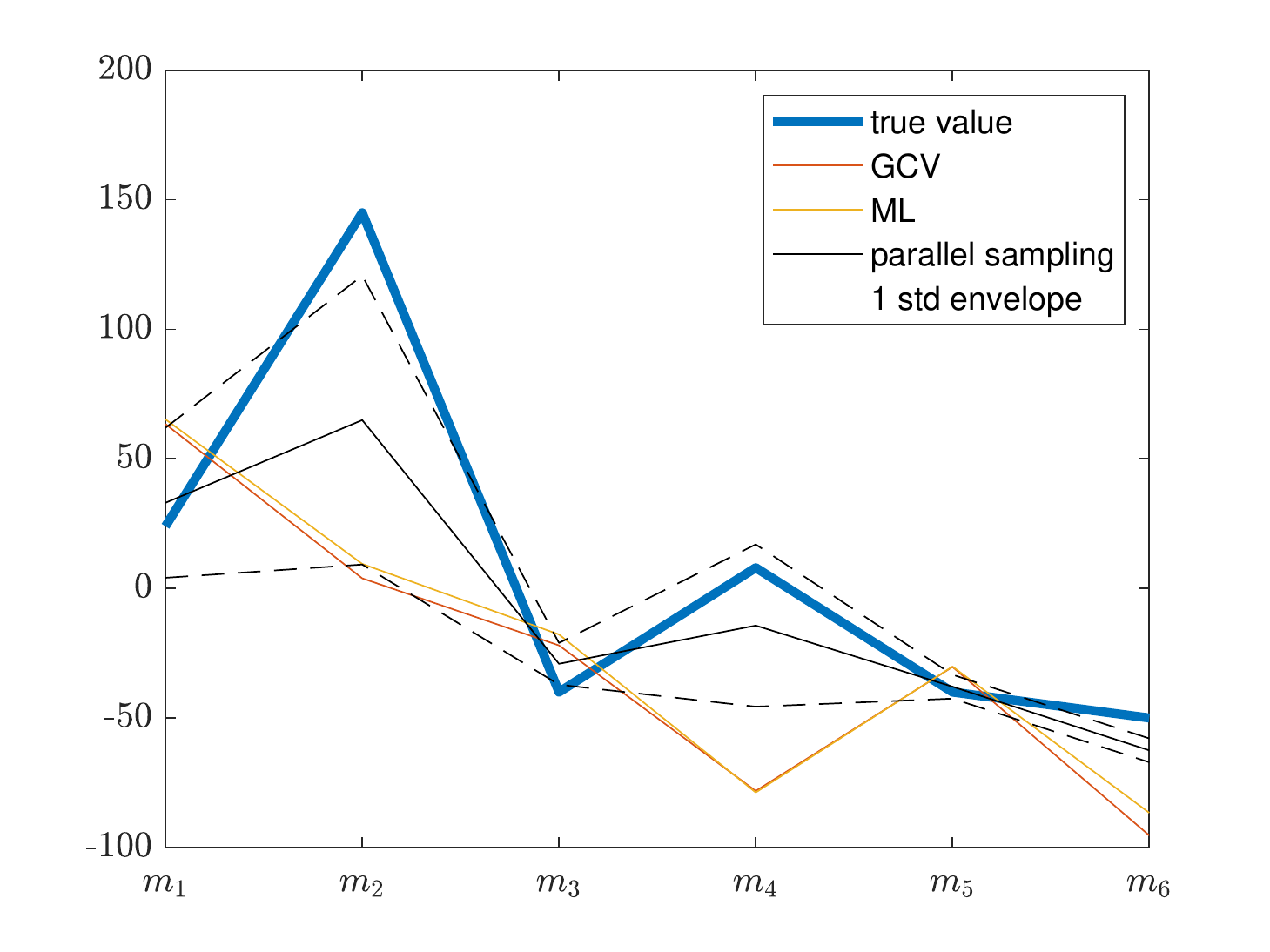} 
			  \includegraphics[scale=.4]{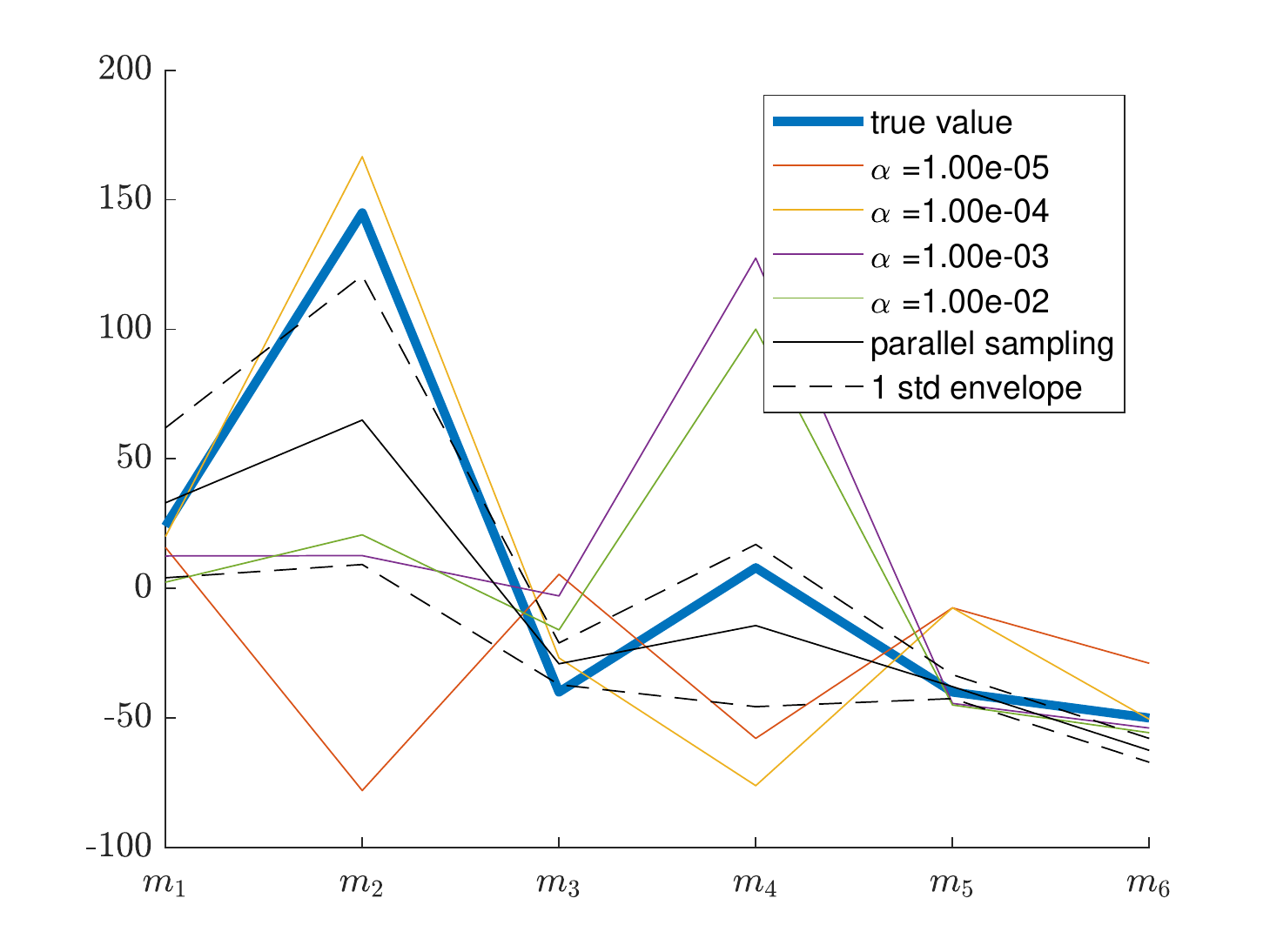} 
				 \caption{Computed values of $\bm$ for the high noise scenario. Top left:
				results obtained by  the GCV  and the ML methods without providing a starting point. Top right: 
				values obtained if these methods
				are started from the favorable point $\tilde{\bm}$. Bottom: 
				the CLS method for different values of $\alpha$ and resulting computed $\bm$.
				In each graph, 
				the true value of $\bm$ and the expected value of $\bm$ with the one standard deviation 
				envelope computed by the parallel sampling method introduced in this paper are indicated for comparison.}
    \label{outcomes high}
\end{figure}

\begin{figure}
   \centering
	  \includegraphics[scale=.4]{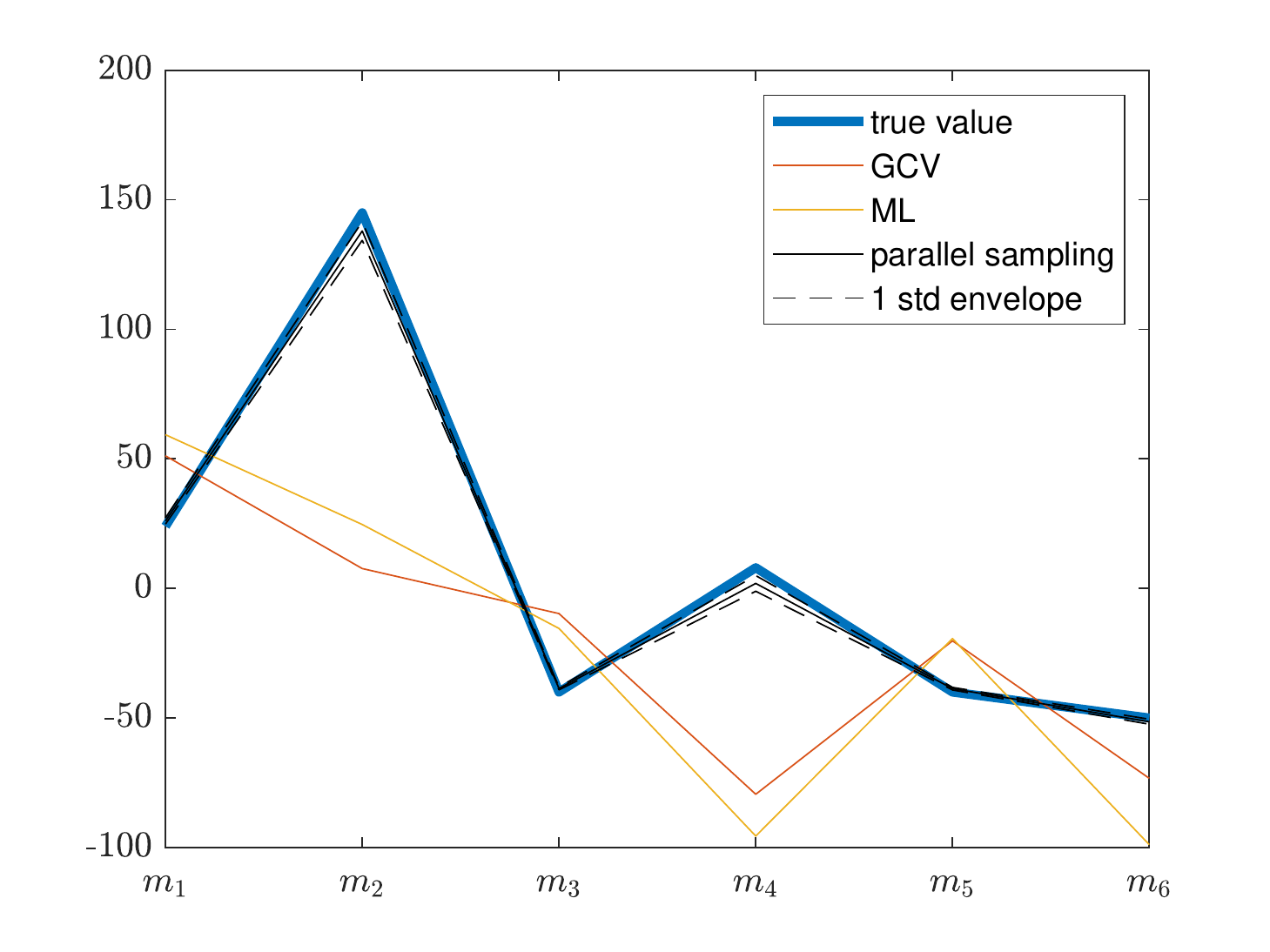} 
			  \includegraphics[scale=.4]{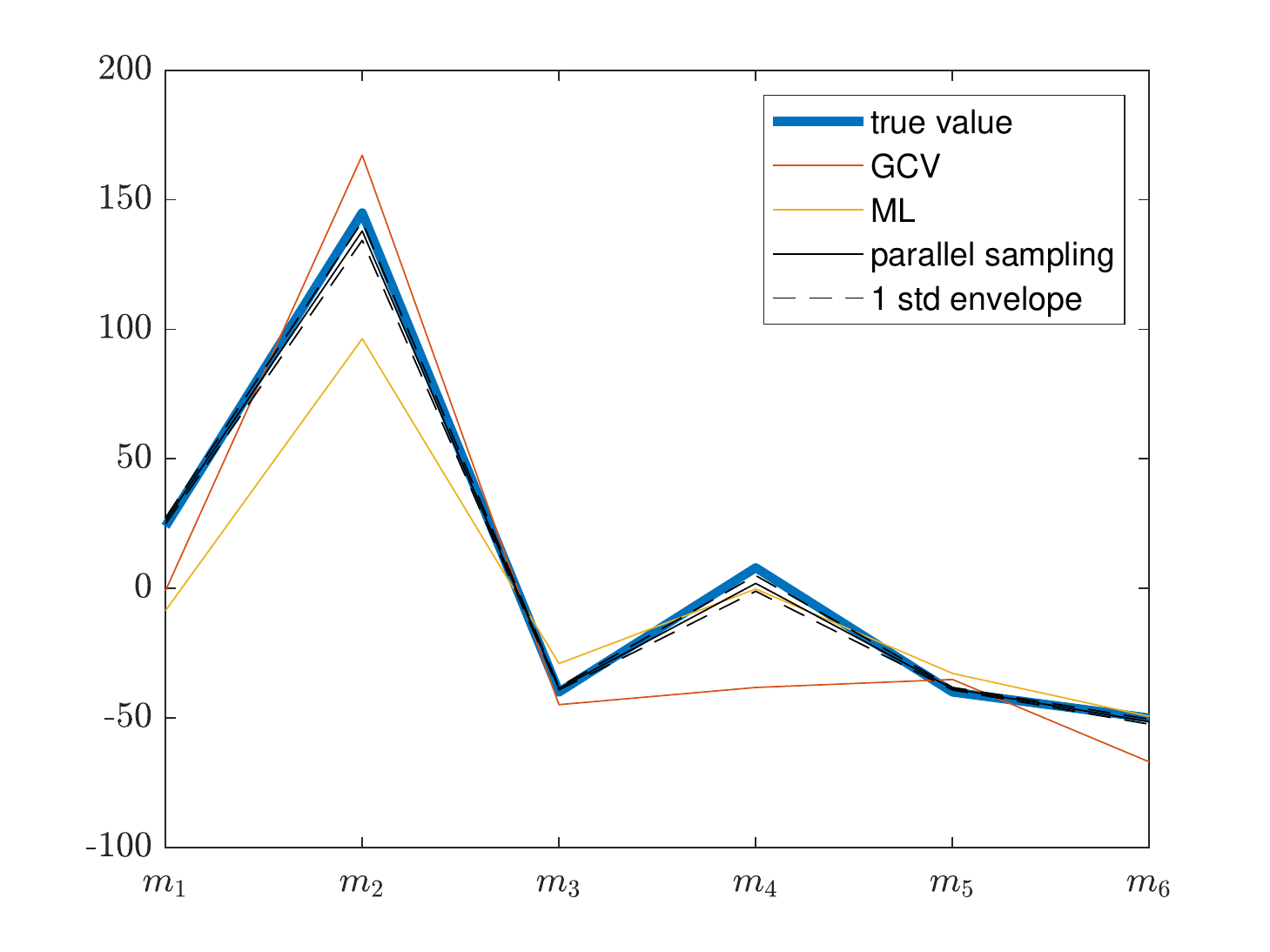} 
			  \includegraphics[scale=.4]{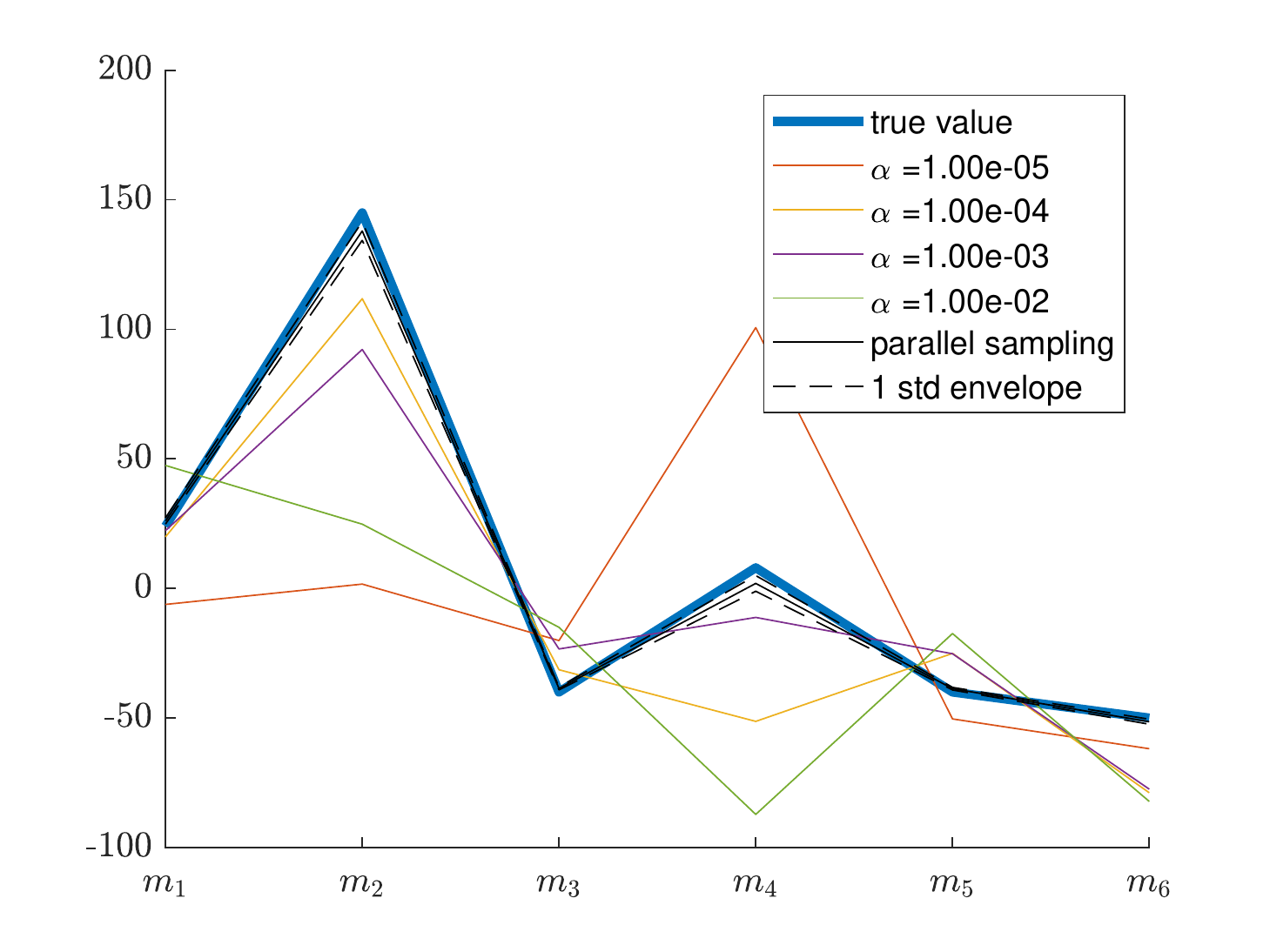} 
			   \caption{Same caption as in Figure \ref{outcomes high} for the low noise scenario. }
    \label{outcomes low}
\end{figure}

For the GCV method, 
we minimized
\bea
 \f{\|  (I_n - A_\bm(A_\bm'A_\bm + \alpha R'R)^{-1} A_\bm' ) \bu \|^2}{ \mbox{tr }(I_n - A_\bm(A_\bm'A_\bm + \alpha R'R)^{-1} A_\bm' )^2},
\eea
and for the ML method,
\bea
 \f{\bu'(I_n - A_\bm(A_\bm'A_\bm + \alpha R'R)^{-1} A_\bm' ) \bu }{(\det (I_n - A_\bm(A_\bm'A_\bm + \alpha R'R)^{-1} A_\bm'  ))^{1/n}}    , 
\eea
as functions in $(\bm, \alpha)$
where 
 $\bm$ is in  the  subset ${\cal B}$ of  $[-200,200]^6$ defined earlier
and $\log_{10} \alpha$ is between -5 and 0.
We used 
the Matlab 2020 function \texttt{surrogateopt}  
to search for the minimum.
This function is based on a minimization algorithm proposed in
\cite{gutmann2001radial} which is specifically designed for problems where function evaluations
are expensive.
This algorithm uses a radial basis function interpolation
to determine  the next point where the objective function should be evaluated.
This is a global non-convex minimization algorithm that uses random restart points in an attempt to 
avoid being trapped in local minima.
We report the computed minima for the high noise scenario in Figure \ref{outcomes high}, top left graph, where we also 
indicate the true value of $\bm$ and the computed expected value 
of $\bm$ obtained 
by our parallel sampling algorithm, together with the one standard deviation 
envelope. We observe that the computed values of $\bm$ obtained by the GCV or the ML algorithm 
are not as close to the true value of $\bm$ as the computed expected value
of $\bm$.
In this high noise scenario, the one standard deviation is  particularly 
informative 
but it cannot be provided by the GCV or the ML algorithm. 
In the second graph of Figure \ref{outcomes high}, we show computed values of $\bm$ if the GCV and the ML algorithms are started from the particularly favorable point
$\tilde{\bm} = (0, 100, -20,  0, -20,  -30)$ and $\tilde{\alpha} = 10^{-4}$.
These two algorithms perform only marginally better despite this favorable head start.
In Figure \ref{outcomes low}, we report results for the low noise scenario. 
It is even clearer in that case that our parallel sampling algorithm outperforms the GCV and the ML 
methods even if they are provided with the favorable starting point
$\tilde{\bm} $.  \\
We also show results obtained by the CLS method in the third graph of Figure 
\ref{outcomes high}
for the high noise scenario, Figure \ref{outcomes low} for the low noise scenario. 
In our particular application $\sigma$  is not known, and even then it would be unclear 
what fixed value to choose for $\alpha$ because of the problem dependence on $\bm$.
It is then common to try   fixing a few values for $\alpha$ 
(we show results for four values of $\alpha$),
and to minimize $\|  A_\bm \bg -  \bu\|^2   + \alpha \| R \bg \|^2$
in $\bm$ and $\bg$. In the high noise scenario, if  $\alpha = 10^{-4}$, this led to
results that are 
better  than 
those obtained by
the GCV or the ML methods. In the low noise scenario,  best results were obtained for 
$\alpha = 10^{-3}$, beating again GCV and ML.
The caveat is that in a real world situation, the  solution to the inverse problem is unknown, so it would be difficult to decide
which  of the four CLS solutions to select.

\section{Conclusion and perspectives for future work}
We have proposed in this paper an algorithm for a finite dimensional 
inverse problem that combines  a  linear unknown $\bg$ 
and a nonlinear unknown $\bm$.
We have used  the maximum likelihood 
assumption for the   prior  of $\bg$ scaled by a parameter $\alpha$
to derive a posterior distribution
for $(\bm, \alpha)$.
Using this posterior distribution we have  built a   parallel  sampling algorithm 
for computing the expected value, the covariance, and 
the  marginal probability distributions of $\bm$.
This algorithm is particularly well suited to the  fault inverse problem where 
$\bm$ models a set of geometry parameters for the fault and 
$\bg$ models a slip on that fault, while the data $\bu$
is sparse and noisy and the linear operator giving $\bu$ from $\bg$ has an unbounded
inverse. Our numerical simulations have shown
that our  parallel  sampling algorithm 
leads to better results than those obtained from minimizing 
the ML, the GCV, or the CLS functionals.
Our algorithm 
  automatically adjusts to a 
	good range for the regularization parameter $\alpha$
	relative   to the  noise level and avoids being trapped in local minima. \\
So far,  our numerical simulations have  focused  on the 
case $q << n << p$, where $\bm$ is in $\RR^q$, the measurements 
$\bu $ are in $\RR^n$, and the forcing term $\bg$ is in $\RR^p$.
However, there are  applications in geophysical sciences where 
measurements are nearly continuous in space and time. This often comes at the price
of higher error margins, so  this would
correspond to the case where $n$ and $p$ are of the same order of magnitude, but $\sigma$ 
is larger, where $\sigma^2 I_n $ is the covariance of the noise. 
Another interesting line of research would be to consider the case where
$q$ is much larger 
which would model an inverse problem
that depends non-linearly on a function, for example the coefficient of a PDE. \\\\
\Large{\bf{Funding}} \\
\normalsize
This work was supported by
 Simons Foundation Collaboration Grant [351025].

\begin{appendices}
 
\section{Proof of formula (\ref{dets})}\label{det proof}
It is well known that the non-zero eigenvalues of $A_\bm'A_\bm$ and $A_\bm A_\bm'$ are the same thus 
$$
\det(\alpha^{-1} A_\bm' A_\bm + I_p) = \det(\alpha^{-1} A_\bm A_\bm' + I_n).
$$
Expanding and simplifying shows that 
\bea
(I_n - A_\bm(\alpha I_p + A_\bm' A_\bm)^{-1}A_\bm')(I_n + \alpha^{-1} A_\bm A_\bm' )=
I_n,
\eea
thus,
\bea
(I_n + \alpha^{-1} A_\bm A_\bm' )^{-1} = I_n - A_\bm(\alpha I_p + A_\bm' A_\bm)^{-1}A_\bm',
\eea
and formula \eqref{dets} is now clear.

\end{appendices}


\end{document}